\documentclass[12pt, a4paper]{article}

\usepackage[utf8]{inputenc}
\usepackage{amsmath, amssymb}
\usepackage{amsfonts}
\usepackage{amsthm}
\usepackage[margin=2cm]{geometry}
\usepackage[colorlinks=true, citecolor=blue, linkcolor=black, urlcolor=black]{hyperref}
\usepackage{graphicx}
\usepackage[dvipsnames]{xcolor}
\usepackage{float}
\usepackage{verbatim}
\usepackage{geometry}
\usepackage{subcaption}
\usepackage{tkz-graph}  
\usepackage{accents}
\usepackage[round]{natbib}
\usetikzlibrary{shapes.geometric}
\usepackage{mathtools}
\mathtoolsset{showonlyrefs}

\newtheorem{theorem}{Theorem}
\newtheorem{lemma}{Lemma}

\newtheorem{proposition}{Proposition}
\theoremstyle{definition}
\newtheorem{example}{Example}[section]

\title{Local asymptotics of selection models with applications in Bayesian selective inference}

\author{
    Daniel García Rasines\textsuperscript{1}\thanks{Address for correspondence. CUNEF Universidad, Calle Almansa 101, Madrid, 28040, Spain. \href{mailto:daniel.garciarasines@cunef.edu}{daniel.garciarasines@cunef.edu}} 
    \and 
    G. Alastair Young\textsuperscript{2}\thanks{Imperial College London}
}

\begin{document}

\maketitle

\abstract{Contemporary focus on selective inference has renewed interest in the theory of selection models. In this paper, we analyze the asymptotic properties
of selection models built on independent and identically distributed observations. We show that, under suitable regularity conditions, they behave asymptotically like a sequence of Gaussian selection models. This provides a natural generalization of the Local Asymptotic Normality framework of \cite{lecam}, and indicates a notion of local asymptotic selective normality as the appropriate simplifying theoretical framework for analysis of selective inference. As a key application, we consider the methodological consequences of the asymptotic theory for Bayesian selective inference. Specifically, we prove that the posterior distribution constructed from a selection model under a fixed prior is asymptotically equivalent to the posterior derived in the corresponding asymptotic Gaussian selection model under a uniform prior. Notably, the latter is often mis-calibrated in a frequentist sense, particularly for one-sided selection mechanisms. This demonstrates that the familiar asymptotic equivalence between Bayesian and frequentist approaches does not hold under selection.}



\maketitle

\section{Introduction}
The techniques of selective inference aim to ensure validity in situations where the same data is used to select the inference to be considered and also to conduct it. Three main approaches have been proposed: the simultaneous approach \citep{berketal}, which ensures validity irrespective of the inference performed and the selection procedure; the `condition on selection' approach \citep{fithianetal}, where inference is conditioned on those hypothetical datasets which might have occurred which would have led to selection of the same inferential problem; and information splitting methods, such as data splitting, where information in the sample is split between that used for selection and that used for inference: see for example \cite{splitting, Leiner, Neufeld}. 

In this work, we study the asymptotic behavior of selection models, which arise as the working inferential models under the conditional and information-splitting approaches to selective inference. These are statistical models in which the data distribution is conditioned on the realization of a certain event. Although selection models have seen renewed attention in recent years, they have long been ubiquitous in statistical inference. Their study dates back at least to \cite{fisher34}, and many of their foundational mathematical properties were developed by \cite{fraser52, fraser66}. These models have traditionally been used to model situations involving sampling bias, where samples are observed only if they satisfy a certain condition---for example, if they belong to a certain portion of the population \citep{rao85}. Selective inference uses a rather more general formulation, which is described in Section \ref{SEC: framework}. 

Within the selective inference literature, much of the methodology for selection models has focused on Gaussian data, typically assuming a known (or consistently estimable) covariance structure. A central reason is that Gaussian models enable a natural treatment of nuisance parameters through standard conditioning arguments; see e.g. \cite{leeetal}. As a result, the theoretical properties of selection models within this setting are well understood, both from the frequentist and the Bayesian positions. The main contribution of this work is the identification of an asymptotic connection between general parametric selection models and Gaussian ones, thereby providing a uniform framework for analyzing the properties of the former and deriving powerful theoretical insights from studying the latter.

More broadly, this paper contributes to the literature on non-standard asymptotics by establishing a limiting behavior which deviates from the standard Gaussian framework. We show that, under mild regularity conditions, a sequence of selection models from a parametric distribution is asymptotically equivalent to a sequence of selection models derived from a Gaussian location model, which, in contrast to standard asymptotics, can be highly asymmetrical if the selection effect is significant (Section \ref{SEC: SLAN}). The relevant Gaussian models inherit a selection mechanism for the mean parameter that often depends on both the original model and the underlying data-generating parameter in an intricate manner, as illustrated in the examples of Sections \ref{SEC: exponential} and \ref{SEC: gaussian}. This showcases the unconventional nature of the analysis, compared to standard asymptotic theory, and underscores the need for careful consideration of the appropriate asymptotic framework for inference. 

Our contribution can be regarded as a generalization of the Local Asymptotic Normality framework of \cite{lecam} to the family of selection models. It indicates a notion of local asymptotic selective normality as the appropriate asymptotic framework for analysis of selective inference problems. In the context of selective inference, therefore, in some substantial generality, a parametric inference problem is asymptotically equivalent to an identifiable selective inference carried out on a location parameter under the Gaussian assumption, conducted from a single observation. As noted above, in contrast to the conventional Local Asymptotic Normality framework, the asymptotic inferential problem considered here can remain rather complicated due to the intricate nature of the relevant asymptotic selection mechanism.

A key motivation for the study of Local Asymptotic Normality is that analysis of the simpler approximating Gaussian model yields information about asymptotic properties in the original parametric model. In this paper, we focus on the methodological implications for Bayesian inference in selection models. As a natural corollary of the asymptotic approximation, we derive an extension of the Bernstein-von Mises Theorem to selection models (Section \ref{SEC: Bayes}). We show that in a selection model with a fixed prior on the parameter, the Bayesian posterior is asymptotically equivalent to that obtained in a Gaussian selection model with a uniform prior on the location parameter. Crucially, this is a non-Gaussian limit, with significant practical consequences for selective inference, as the posterior in such a Gaussian selection model is mis-calibrated from a frequentist perspective. In these models, frequentist calibration requires specification of a prior density which can depend strongly on the sample size, as discussed in \cite{handbook}.

\section{Framework}\label{SEC: framework}

The term selection model may refer to slightly different types of statistical models, depending on the context. For the purposes of this work, a selection model is defined by two elements: a \textit{base model}, understood as the pre-selection model and assumed to contain the true data-generating distribution; and a \textit{conditioning random variable}, upon which the data is conditioned to obtain the selection model. This variable
encodes the \textit{selection mechanism}, i.e. the rule determining which realizations from the base model are observed, or selected, for some purpose. Operationally, a sample $y^n$ is retained only if a certain value of the conditioning variable is observed. This framework is sufficiently general to encompass all applications in the conditional and information-splitting approaches to selective inference, as illustrated in Section \ref{SEC: SLAN}. 

Typically, the base model is well understood from an inferential standpoint, and the objective is to extrapolate its properties to the selection model by analyzing how the selection mechanism interacts with it. In this paper, we will restrict attention to base models consisting of random samples from low-dimensional, regular parametric families. Other classes of base models of central importance in selective inference which are not covered in this work include fixed-design regression, high-dimensional, and non-parametric models.

\subsection{Base model}

Let $Y^n = [Y_1, \ldots, Y_n]$ be a vector of independent and identically distributed (IID) random variables from a distribution $F_\theta$ on a measurable space $(\mathcal{Y}, \mathcal{A})$, where $\theta$ is a parameter from an open set $\Theta \subseteq \mathbb{R}^p$. The model for the full sample, defined over $(\mathcal{Y}^n, \mathcal{A}^n)$, is $\mathcal{M}_n = \{F_\theta^n \colon \theta\in \Theta\}$. We regard $\{\mathcal{M}_n \colon n\geq 1\}$ as the sequence of base models, from which a sequence of selection models will be obtained. It is well known that, under mild regularity conditions, $\mathcal{M}_n$ can be asymptotically approximated by the Gaussian location model $\{N(h, I_\theta^{-1})\colon h\in \mathbb{R}^p\}$, where $I_\theta$ is the Fisher information for a single sample. Specifically, if the model $\{F_\theta : \theta \in \Theta\}$ is differentiable in quadratic mean (DQM; \citealp[][p.~93]{vandervaart}), then the sequence $\{\mathcal{M}_n : n \ge 1\}$ is Locally Asymptotically Normal (LAN) in the sense of \citet{lecam}. That is, the sequence of local log-likelihood ratios admits the expansion
\begin{equation*}
    \log \prod_{i=1}^n \frac{f_{\theta + h_n/\sqrt{n}}(Y_i)}{f_\theta(Y_i)}
    = h^\top I_\theta Z_n - \tfrac{1}{2} h^\top I_\theta h + o_{F_\theta}(1) \quad \text{as } n\to\infty,
\end{equation*}
for every converging sequence $h_n\to h$, where $Z_n \rightsquigarrow N(0, I_\theta^{-1})$ under $F_\theta$, and $f_\theta(y)$ is the density of $F_\theta$ with respect to some fixed measure $\mu$. We use $\rightsquigarrow$ to denote convergence in distribution. In the limit, the local likelihood ratios are distributed as the likelihood ratios of $\{N(h, I_\theta^{-1})\colon h\in \mathbb{R}^p\}$. This translates into an asymptotic equivalence between both models in likelihood-based inference. Many key asymptotic results can be derived from the LAN property, including the asymptotic distribution of maximum likelihood estimators and the Bernstein–von Mises theorem in Bayesian inference; see \cite{vandervaart} for a detailed account. Throughout this work, we shall assume that $\{F_\theta : \theta \in \Theta\}$ is differentiable in quadratic mean. This ensures that the sequence of base models is LAN, which forms the foundation of the asymptotic theory for selection models.

\begin{example}\label{EX: exp}    
    We will illustrate the framework through an exponential base model, where $Y_i \sim \text{Exp}(\theta)$ and $\Theta = (0, \infty)$. Upon setting $Z_n = \theta^2 \sqrt{n}(1/\theta - \bar{Y}_n)$ to be the normalized score at $\theta$, where $\bar{Y}_n$ is the sample mean, the LAN expansion of $\mathcal{M}_n$ specializes to 
\begin{equation*}
    \log \prod_{i=1}^n \frac{f_{\theta + h_n/\sqrt{n}}(Y_i)}{f_\theta(Y_i)}
    = \sqrt{n} h \left(\frac{1}{\theta} - \bar{Y}_n\right) - \frac{h^2}{2\theta^2} + o_{F_\theta}(1) \quad \text{as } n\to\infty.
\end{equation*}
\end{example}

\subsection{Selection model}

A selection model can be obtained from a base model $\mathcal{M}_n$ by conditioning the data on a random quantity $U_n$ such that $U_n\mid Y^n = y^n$ is distribution-constant for every $y^n$. This implies that $U_n$ is either a statistic, such as $u(Y^n) = \mathbf{1}(\bar Y_n > 0)$, where $\mathbf{1}(A)$ denotes the indicator function of the event $A$, or a randomized function of the data, $u(Y^n, W_n)$, say, where $W_n$ is user-generated noise not depending on $\theta$. Although this might seem somewhat artificial, randomizing the selection step is standard practice in selective inference. The purpose is to increase the information about $\theta$ contained in the selection model, in a way analogous to data splitting. For a realized value $U_n = u_n$, the conditional density of $Y^n\mid u_n$ is 
\begin{equation}
    f_{\theta}(y^n\mid u_n) = \frac{f(u_n\mid y^n)}{f_\theta(u_n)} \prod_{i = 1}^n f_\theta(y_i),
\end{equation}
where $f(u_n\mid y^n)$ and $f_\theta(u_n)$ are the conditional and marginal densities of $U_n$ relative to some dominating measures, and $f_\theta(u_n)$ is assumed to be positive for every value of $\theta$. We denote the corresponding probability distribution by $F_\theta^{u_n}$, and define the selection model with base model $\mathcal{M}_n$, conditioning variable $U_n$, and realized $U_n = u_n$ as $\mathcal{M}_n(u_n) = \{F_n^{u_n}\colon \theta\in\Theta\}$. Furthermore, we define a \textit{Gaussian selection model} as any selection model whose base model is a Gaussian location model: $\mathcal{M}_n = \{N(\theta, \Sigma)^n\colon \theta\in \mathbb{R}^p\}$, where $\Sigma$ is known and positive definite. We will show that Gaussian selection models arise as asymptotic first-order approximations of general selection models under suitable regularity conditions.

Our objective is to understand the asymptotic behavior of $\{\mathcal{M}_n(u_n)\colon n\geq 0\}$ for fixed sequences of $U_n$-realizations $\{u_n\colon n\geq 1\}$. For this reason, we will occasionally suppress $u_n$ from the notation. In particular, we will write
\begin{equation}
    f_{\theta}(y^n\mid u_n) = \frac{p_n(y^n)}{\varphi_n(\theta)} \prod_{i = 1}^n f_\theta(y_i),
\end{equation}
where $p_n(y^n) = f(u_n\mid y^n)$ and $\varphi_n(\theta) = f_\theta(u_n)$. This simplifies the expressions, in particular in Section \ref{SEC: Bayes} where we discuss Bayesian selective inference, and aligns better with the notation conventionally employed in most selective inference applications. Henceforth, we will refer to $p_n(y^n)$ as the \textit{selection function}. This encodes all relevant distributional information about the selection mechanism, in the sense that any two statistics and realized values $U_{n1} = u_{n1}$ and $U_{n2} = u_{n2}$ with proportional selection functions $p_{n1}(y^n) \propto p_{n2}(y^n)$ induce the same selection model on $Y^n$. In many applications, $U_n$ is a discrete random variable which indicates whether a certain event of interest $A$ has occurred: $U_n = \mathbf{1}(A)$. In such cases, $A$ is called the \textit{selection event} and, for $u_n = 1$, $\varphi_n(\theta) = P_\theta(A)$ is the \textit{selection probability}. Note, however, that whenever the selection function is bounded, the same selection model can be represented as one induced by a discrete selection variable. Indeed, we may normalize $p_n$ to obtain a genuine conditional probability by defining the equivalent selection function
\begin{equation*}
    \tilde p_n(y^n) = \frac{p_n(y^n)}{\Vert p_n \Vert_\infty}, \quad \Vert p_n \Vert_\infty = \sup\{p_n(y^n)\colon y^n\in \mathcal{Y}^n\},
\end{equation*}
so that $0 \le \tilde p_n(y^n) \le 1$ for all $y^n$. We can then introduce an auxiliary binary variable by letting $\tilde U_n\mid y_n \sim \text{Bernoulli}\{\tilde p_n(y^n)\}$, so that conditioning on $\tilde U_n = 1$ yields a selective model for $Y^n$ that is equivalent to the original one.

Selection models play a central role in the study of sampling bias and selective inference. In the context of sampling bias, it is often assumed that the selection function factorizes into marginal selection functions for each observation, i.e. that $p_n(y^n) = \prod_{i = 1}^n p(y_i)$, in which case the resulting selection model is IID and standard asymptotics apply under mild conditions. In selective inference, on the other hand, the selection condition typically acts on the dataset as a whole, and such a representation does not hold. Under the conditional approach to selective inference, if a subparameter $\psi = g(\theta)$ is only analyzed provided a pre-specified condition on the sample is satisfied, then inference on $\psi$ ought to be carried out conditionally on this event, commonly referred to as the \textit{selection event}. Here, pre-specified means that the condition was decided by the statistician before analyzing the data. In the notation introduced above, the selection event takes the form $\{U_n = u_n\}$ for a specific value $u_n$, and its occurrence typically indicates that $\psi$ is worth investigating for a particular reason. For example, this could be rejection of some null hypothesis, indicating a significant effect size. We will see some common examples of selection events in the following section. It is worth noting that the conditioning operation in this case is not intrinsic to the data sampling mechanism, but a formal distributional correction undertaken to restore repeated-sampling validity to the inferences. 

\setcounter{example}{0} 
\begin{example}[ctd.]
Consider the exponential base model of Example \ref{EX: exp} and assume that, for scientific reasons, we are only interested in $\theta$ if it is large. That is, we would like to perform inference on $\theta$ only if the sample information indicates a large enough rate $\theta$, according to some measure of calibration. Let $U_n = \mathbf{1}(\sqrt{n}\bar{Y}_n + W < t)$ for some threshold $t$, where $W\sim N(0, 1)$ independently of $Y^n$, and assume a sample is selected for further analysis only if $U_n = 1$, which may be regarded as rejecting the null hypothesis in a randomized one-sided test for $\theta$. The conditional approach to selective inference dictates that inference following observation of $U_n = 1$ ought to be based on the corresponding selection model, with respective selection function and probability
\begin{eqnarray}
    p_n(y^n) &=& \Phi\left( t - \sqrt{n}\bar y_n   \right); \\
    \varphi_n(\theta) &=& E_\theta\left[ \Phi\left( t - \sqrt{n}\bar Y_n   \right) \right];
\end{eqnarray}
where $\Phi$ is the $N(0, 1)$ distribution function.

Suppose we wish to test $H_0\colon \theta\leq \theta_0$ against the alternative $H_1 \colon \theta > \theta_0$ for large $n$. An exact test under the selection model requires comparison of the observed $\bar y_n$ with the sampling distribution of $\bar Y_n$ when $\theta = \theta_0$, given selection. Analytically, this requires consideration of the convolution of a Gamma and a Gaussian distribution, which in this case is analytically complicated, though easily simulated. In the following section we will show local asymptotic equivalence, at a given $\theta$, between $\mathcal{M}_n(1)$ and a Gaussian selection model. The test of $H_0$ can therefore be approximated by referring the observed value of $z_n(\theta_0)=\theta_0^2 \sqrt{n}(1/\theta_0 - \bar{y}_n)$ to the asymptotic density
\begin{equation} \label{asydensity} f_{\theta_0}(z)= \frac{1}{\theta_0} \cdot \frac{\phi(\frac{z}{\theta_0})\Phi \bigl(t-\sqrt{n}/\theta_0+\frac{z}{\theta_0^2}\bigr)}{\Phi\big\{(t-\sqrt{n}/\theta_0)/\sqrt{1+1/\theta_0^2}\bigr\}}, \quad z \in \mathbb{R},
\end{equation}
where $\phi$ is the standard Gaussian density. Evidence against $H_0$ is provided by small values of $z_n(\theta_0)$. The appropriate approximate $p-$value is now easily evaluated by a simple numerical integration of the density \eqref{asydensity}. 

Such use of the asymptotic Gaussian selection model to examine repeated sampling properties in the initial exponential selection model is illustrated in  Figure \ref{pvalueplots} for two sample sizes, $n = 50, 200$. For each sample size, we consider, over a series of 10,000 replications, the distribution in the two models of the $p-$value for testing $H_0\colon \theta\leq 2$ against the alternative $H_1 \colon \theta > 2$, under two scenarios, when $H_0$ is true, and when the alternative $H_1$ is true, with $\theta = 2.15$. In each case we calculate the $p-$value under the true selection model by simulating the null sampling distribution of $\bar Y_n$ given the selection condition $\sqrt{n}\bar{Y}_n + W < t$ from a series of 100,000 samples. The asymptotic $p-$values are obtained by numerical integration of the density \eqref{asydensity}. The selection thresholds are set as $t=2.1$ and $t=5.6$ respectively for the two sample sizes $n=50, 200$: these thresholds correspond to low selection probability $ \approx 10\%$ under the null. Figure \ref{pvalueplots} shows that even for $n=50$ the two $p-$value distributions are close under both the null and alternative. The distributions are quite indistinguishable for the larger sample size: repeated sampling behaviour under the original randomized selection model is closely matched by that seen in the simpler asymptotic selection model. Inferential properties under the true selection model are accurately predicted from those of the asymptotic approximating Gaussian selection model.

\begin{figure}[h]
    \centering
    \includegraphics[width=\textwidth]{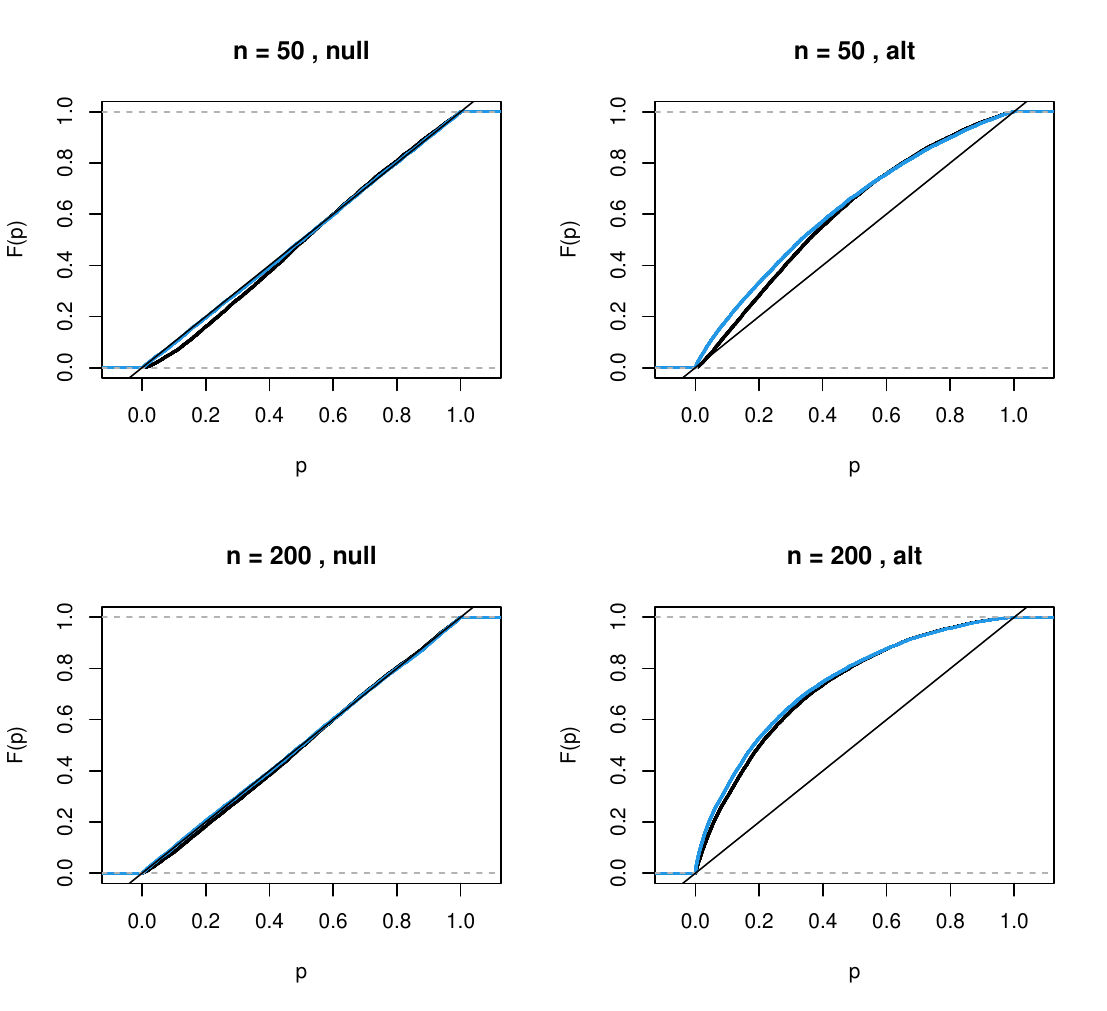}
    \caption{Repeated sampling distribution of $p-$values under the true randomized exponential selection model (blue) and the asymptotic Gaussian selection model (black), $n=50,200$, under null and alternative.}
    \label{pvalueplots}
\end{figure}


In the Bayesian setting, we will show that the repeated-sampling behavior of posterior distribution for $\theta$ in the selection model can be derived from the asymptotically equivalent Gaussian model. This reveals a fundamental divergence between the frequentist and Bayesian approaches in selection models that contrasts with classical theory for IID models.

\end{example}

\section{Asymptotic expansion of the likelihood} \label{SEC: SLAN}

In this section, we establish a formal asymptotic connection between the sequence of selection models, $\{\mathcal{M}_n(u_n)\colon n\geq 1\}$, and a sequence of Gaussian selection models. At the heart of the result is an asymptotic selection function that reflects how selection influences the model asymptotically. This can be regarded as providing an extension of the standard LAN framework, which can be recovered by considering a constant selection function, $p_n(y^n)\propto 1$. 

\subsection{Main result}
First, we examine how selection acts on Gaussian models. Let $Z\sim N(h, \Sigma)$, with parameter $h\in \mathbb{R}^p$ and fixed positive definite covariance $\Sigma \in \mathbb{R}^{p\times p}$, and consider the selection model induced by an arbitrary selection function $p(z)$. The selective log-likelihood ratio for this model around $h = 0$ is given by
\begin{equation} \label{EQ: exp_Gauss}
\log \frac{f_{h}(Z)}{f_{0}(Z)} = h^\top \Sigma^{-1} Z - \frac{1}{2} h^\top \Sigma^{-1}h   + \log \frac{\varphi(0)}{\varphi(h)}, \quad \varphi(h) = E_h[p(Z)].
\end{equation}
Note that the selective log-likelihood ratio follows itself a selective Gaussian distribution under any $h\in \mathbb{R}^p$, just as the log-likelihood ratio of a conventional, non-selective Gaussian model is itself Gaussian. 

The main result of this paper, Theorem \ref{SLAN}, establishes conditions under which the sequence of selective log-likelihood ratios of a given selection model $\mathcal{M}_n(u_n)$ admits an asymptotic expansion of the form \eqref{EQ: exp_Gauss} for a certain specification of selection function, which in general depends on the true generating parameter $\theta$. It is worth remarking that, since selection functions can vary arbitrarily with $n$, $\{\mathcal{M}_n(u_n)\colon n\in \mathbb{N}\}$ does not in general converge to a fixed Gaussian selection model. Instead, it can, under suitable conditions, be approximated by a sequence of such. 

Before stating the Theorem, we need to introduce some notation.

\textbf{Notation}.  In the non-selective model $\mathcal{M}_n$, let $l_\theta(y_i) = \log f_\theta(y_i)$ denote the log-likelihood for the $i$-th observation, $\triangledown l_\theta(y_i)$ the score, and $I_\theta = E_\theta [\triangledown l_\theta(Y_i)\triangledown l_\theta(Y_i)^\top]$ the per-observation Fisher information, assumed to be positive definite for all $\theta\in \Theta$. For a function $g: I \subseteq \mathbb{R}^p \to \mathbb{R}$, consider its Lipschitz norm,
\begin{equation}
\Vert g\Vert_\text{L} = \sup \left\{ \frac{\vert g(x) - g(y) \vert}{\Vert x - y \Vert} \colon x, y\in I, x\neq y \right\},
\end{equation}
where $\Vert \cdot \Vert$ is the Euclidean norm; its supremum norm, $\Vert g \Vert_\infty = \sup \{ \vert g(x)\vert  \colon x\in I\}$; and its bounded-Lipschitz norm, $\Vert g\Vert_{\text{BL}} = \Vert g \Vert_\infty + \Vert g\Vert_\text{L}$. For a $J\subseteq I$, let $g_{|J}$ be the restriction of $g$ to $J$. The bounded-Lipschitz distance between the distributions of two $p$-dimensional random vectors $X$ and $Y$ is
\begin{equation}
     \Vert X - Y \Vert_{\text{BL}}= \sup \left\{ E_X[g(X)] - E_Y[g(Y)] \colon\quad g\colon \mathbb{R}^p\to \mathbb{R},\quad   \Vert g\Vert_{\text{BL}}  \leq 1 \right\}.
\end{equation}
For a probability distribution $\mu$ on $\mathbb{R}^p,$ consider its Lebesgue decomposition: $\mu = \mu_{C} + \mu_S$, where $\mu_C$ is absolutely continuous with respect to the Lebesgue measure and $\mu_S$ is concentrated on a set of null Lebesgue measure. We say that $\mu$ has an absolutely continuous component if $\mu_S(\mathbb{R}^p) < 1$. Finally, the Total Variation distance between two probability distributions $\mu$ and $\nu$ on $\mathbb{R}^p$ is $\Vert \mu - \nu \Vert_{\text{TV}} = \sup\{\vert \mu(A) - \nu(A) \vert\colon A\in \mathcal{B}(\mathbb{R}^p)\}  $, where $\mathcal{B}(\mathbb{R}^p)$ is the Borel sigma algebra on $\mathbb{R}^p$.

We are now equipped to state the main result of the paper.

\begin{theorem} \label{SLAN}
    Suppose that $\{F_\theta\colon \theta\in \Theta\}$ is DQM, define $Z_n \equiv Z_n(\theta) = n^{-1/2} I_\theta^{-1}  \sum_{i = 1}^n \triangledown l_\theta(Y_i) $, and assume that the sequence of selection functions is uniformly bounded:
\begin{equation}
\sup\{\Vert p_n \Vert_\infty\colon n\geq 1\} <\infty.
\end{equation}    
For a fixed $\theta\in \Theta$, assume that $M(x) = E_\theta[\exp\{ x^\top \triangledown l_\theta(Y_1)\}] < \infty$ in a neighborhood of the origin, that there exists a sequence of measurable functions $\{p^*_n(\cdot ; \theta)\colon \mathbb{R}^p \to [0, \infty) \colon n\in \mathbb{N}\}$ satisfying $p^*_n(z; \theta) = E_\theta \left[p_n(Y^n) \mid  Z_n = z  \right]$ for all $z$ in the support of $Z_n$, and that either of the following conditions are met:
       \begin{enumerate}
        \item \label{C1} The distribution of $\triangledown l_\theta(Y_1)$ has an absolutely continuous component;
        \item \label{C2} $\sup\{\Vert p^*_n(\cdot ; \theta) \Vert_{\mathrm{BL}}\colon n\geq 1\} <\infty$.
    \end{enumerate} 
Define
\begin{equation}
    \varphi^*_n(h) = E[ p^*_n(Z; \theta)], \quad Z\sim N(h, I_\theta^{-1}).
\end{equation} 
If $\inf \{\varphi_n(\theta)\colon n\geq 1\} > 0$, then, for all convergent sequences $h_n\to h$,
\begin{equation}
     \log \frac{f_{\theta + h_n/\sqrt{n}}(Y^n\mid u_n)}{f_{\theta}(Y^n\mid u_n)} = h^\top I_\theta Z_n - \frac{1}{2} h^\top I_\theta h + \log\frac{\varphi_n^*(0)}{\varphi_n^*(h)} + o_{F^{u_n}_\theta}(1).
\end{equation}
Furthermore, if condition 1 is met, $\Vert Z_n - Z_n^* \Vert_{\mathrm{TV}} \to 0$ conditionally on $u_n$, where $Z_n^*$ follows a selective $N(0, I_\theta^{-1})$ distribution with selection function $p^*_n(z; \theta)$. Otherwise, $\Vert Z_n - Z_n^* \Vert_{\mathrm{BL}} \to 0$.
    
\end{theorem}

An asymptotic connection is therefore established between the original model and a Gaussian selection model with base distribution $N(h, I_{\theta}^{-1})$ and selection function $p_n^*(z;\theta)$. In particular, the result implies that, to first order, the selective likelihood depends on the data only through the normalized score $Z_n$, so $Z_n$ is asymptotically sufficient and provides a natural asymptotic pivot for inference. Testing applications are explored in later sections.

The proof of Theorem 1 can be found in the supplementary material along with the proofs of all the other theoretical results of the paper. It takes as starting point the LAN property of the underlying, non-selective model, from which the first two terms of the expansion can be derived. The main technical difficulty lies in validating the approximation of the probability ratios $\varphi_n(\theta + n^{-1/2}h)/\varphi_n(\theta)$ by their Gaussian counterparts. This requires consideration of the moment generating function of $Z_n$ given selection. The existence of $M(x)$, together with the regularity assumptions on the selection functions in the limiting Gaussian models, prove to be sufficient to ensure uniform validity. While the regularity conditions on the selection functions are relatively mild, they do exclude scenarios where the probability of selection vanishes asymptotically. In such cases, the Gaussian approximation may fail to capture the behavior of the posterior or likelihood ratios, and the theoretical guarantees of Theorem 1 no longer apply. We therefore add a note of caution that the current framework is restricted to selection events with probabilities bounded away from zero.

It is also important to note that the selection function in the asymptotic model, $p^*_n(z; \theta)$, depends on the true parameter $\theta$. This dependency introduces additional complexity to the asymptotic analysis. Under the standard LAN framework, the true value of $\theta$ only influences the covariance of the limiting Gaussian location model, leaving the group structure unaltered and facilitating a unified approach to analyzing the properties of inferential procedures. By contrast, under selection, $\theta$ has a deeper structural influence on the limiting model through its presence in $p^*_n$. This may affect the asymptotic distribution of estimators and test functions, as well as optimal choice of inferential approach. It remains to be explored whether the true asymptotic selection function can be effectively approximated by $\hat p^*_n(z) = p^*_n(z; \hat \theta)$, where $\hat\theta$ is a $\sqrt{n}$-consistent estimator of $\theta$.

\subsection{Vanishing selection probabilities}

The requirement that $\varphi_n(\theta)$ is bounded away from zero is needed to ensure that the Gaussian approximation of the normalized score remains valid after conditioning on selection and to preserve the scale of likelihood ratios. This condition guarantees asymptotic regularity of the model even if selection is deterministic, that is, when the selection functions are of the form $p_n(y^n) = \mathbf{1}(y^n\in A_n)$ for some events $A_n$. However, it is well known that, in these scenarios, if $y^n$ falls close to the boundary of $A_n$, irregularities can arise, leading to poor inferential performance unless $n$ is very large.

Scenarios with $\varphi_n(\theta)\to 0$ can arise in practice, particularly in models involving many parameters. A clear example arises when inferring the mean of a Gaussian model when the true mean is negative and the selection region is $A = (t_n, \infty)$ for some sequence $t_n\geq 0$. Extending the framework to handle such vanishing-probability cases is an important direction for future work. While this asymptotic regime is not covered by the Theorem, in some scenarios it is possible to verify the boundedness condition either exactly or asymptotically:
\begin{itemize}
    \item[(i)] If the selection condition is rejection of a simple hypothesis $H_0\colon \theta = \theta_0$ at a level $\alpha$, then $\varphi_n(\theta) \geq \alpha$ for all $n$, provided the test is unbiased. Thus, as long as $\alpha$ does vanish as $n\to\infty$, the condition is satisfied. This is typically the case for multiple comparisons models, such as the one described in Section \ref{SEC: gaussian}, with a bounded number of parameters.       
    \item[(ii)] If the selection mechanism  allows for the existence of a sequence of $\sqrt{n}$-consistent estimator $\hat\theta_n$ of $\theta$ conditionally on selection, then, for large values of $n$, the magnitude of $\varphi_n(\theta)$ can be assessed through that of $\varphi_n(\hat\theta_n)$ (Lemma \ref{lem}). This is always possible if selection acts on a subset of the observations, leaving the remaining ones free of selection bias (see Example \ref{E2} below), and often possible if it acts on a randomized version of the dataset (Examples \ref{E3} and \ref{E4}). Examples of randomization mechanisms that enable estimators with the required accuracy are provided by \cite{tian}, \cite{Leiner}, \cite{Neufeld}, and \cite{Dharamshi}.
    
\end{itemize}

\begin{lemma}\label{lem}
    If $\{F_\theta\colon \theta\in \Theta\}$ is DQM and there exists a sequence of estimators $\hat\theta_n$ such that $n^{1/2}\Vert \hat\theta_n - \theta \Vert = O_{F^{u_n}_\theta}(1)$, then $\varphi_n(\theta) = o(1) \iff \varphi_n(\hat\theta_n) = o_{F^{u_n}_\theta}(1)$.  
\end{lemma}

Selection rules that do not allow for the existence of $\sqrt{n}$-consistent estimators are typically deterministic and yield inferences with very limited power. Their lack of power often requires the collection of additional data to draw meaningful conclusions about the parameter, leading to a data carving scenario (Example \ref{E2}) in which estimation of $\theta$ with the required accuracy is feasible.

\subsection{Examples of common selection mechanisms}

The following examples are representative of the main selection mechanisms encountered in selective inference.

\begin{example}\label{E1}
\textbf{Deterministic selection}.
Let $Y_1, \ldots, Y_n$ be an IID sample from an exponential family distribution with density 
\begin{equation}
    f_\theta(y_i) = h(y_i)\exp\left\{ 
 \theta^\top y_i - A(\theta)\right\}, \quad  \theta\in \Theta.
\end{equation}
Let $U_n = \mathbf{1}(\bar{Y}_n \in E_n)$, where $\bar{Y}_n = n^{-1}(Y_1 + \ldots + Y_n)$ and $E_n\subseteq \mathbb{R}^p$ is an arbitrary sequence of Borel sets with a bounded number of connected components. For example, if $p = 1$ and $E_n = [t_n, \infty)$, such a selection condition corresponds to rejection of a null hypothesis about $\theta$ under a uniformly most powerful test. In this case, $Z_n = \sqrt{n} \triangledown^2 A(\theta)^{-1}\{\bar{Y}_n - \triangledown A(\theta)\}$, and
$p^*_n(z; \theta) =  \mathbf{1}\{n^{-1/2} \triangledown^2 A(\theta) z + \triangledown A(\theta) \in E_n\}$. 
\end{example}

\begin{example}\label{E2}
\textbf{Data carving}.   
Consider the same setting as in the previous example, but with selection condition given by $\bar{Y}^{(1)}_n \in E_n$, where $\bar{Y}^{(1)}_n$ is the mean of the first $n/2$ samples, where $n$ is assumed to be even. Then, 
\begin{equation}
   p^*_n(z; \theta) = \frac{1}{f_{ Z_n}(z)} \int_{E_n} f_{\bar{Y}^{(1)}_n} \{2[n^{-1/2} \triangledown^2 A(\theta) z + \triangledown A(\theta)] - \bar y_1\}  f_{\bar{Y}^{(1)}_n}( \bar y_1 )\mathrm{d}\bar y_1.
\end{equation}
Condition 2 of the Theorem is satisfied provided $\sup\{\Vert \log f_{\bar{Y}^{(1)}_n} \Vert_{\text{L}}\colon n\geq 1\} < \infty$. Conditioning on selection based on a subsample of the data was introduced by \cite{fithianetal} and is commonly referred to as data carving. The main operational consequence of data carving is a smoothing of the selection function.
\end{example}

\begin{example}\label{E3}
\textbf{Randomization}.
Consider the same setting, but suppose now that the selection mechanism acts on $U_n = \sqrt{n} \bar{Y}_n + W$, where the noise $W$ has a known density $f_W(w)$, so that selection takes the form $U_n\in E_n$. Then, $p^*_n(z; \theta) =  P\{\sqrt{n} \bar{y}_n(z) + W \in E_n\}$, so
\begin{equation}
\triangledown p^*_n(z; \theta) = \int_{E_n} -\triangledown^2  A(\theta) \triangledown f_W\{u - \triangledown^2 A(\theta) z - \sqrt{n} \triangledown A(\theta)\} \mathrm{d}u,
\end{equation}
which is bounded provided $f_W(w)$ is bounded. Proposed by \cite{tian}, randomized selection has become common in practice in post-selection inference, as it enables more powerful analyses. Similarly to data carving, the randomization here tends to smooth the selection function
\end{example}

\begin{example}\label{E4}
\textbf{Condition on randomized statistic}.
Consider the same setting as in the previous example. In cases where knowledge of the precise form of $E_n$ is unavailable to the statistician, or too complex to be conditioned on, it is common to condition on the observed value $U_n$ instead \citep{splitting, Leiner, Neufeld}. In that case, the selection function becomes $p^*_n(z; \theta) =  f_W\{u_n - \triangledown^2 A(\theta) z + \sqrt{n} \triangledown A(\theta)\}$, which satisfies the required regularity conditions for most standard noise distributions, such as the standard Gaussian or Laplace distributions employed by \cite{tian}.
\end{example}

\subsection{Exponential model (ctd.)} \label{SEC: exponential}

 We will illustrate the four previous types of selection mechanism in the simple model of Example \ref{EX: exp}. A deterministic selection mechanism would be of the form $\bar Y_n < t_n$ for some predetermined threshold value $t_n$, which might have been obtained through a power analysis or in some other data-independent way. For a given $\theta$ and realization $\bar{Y}_n = \bar{y}_n$, the normalized score for this model is $z = \theta^2 \sqrt{n}(1/\theta - \bar{y}_n)$, and the selection function in the asymptotic model becomes $p^*_n(z; \theta) = \mathbf{1}\{z > \theta^2 \sqrt{n}(1/\theta - t_n)\}$. If the analogous selection rule is applied only to half the samples, so that inference is performed whenever $\bar{Y}_n^{(1)} < t_n$ (data carving), then, assuming $n$ is even and noticing that $Z_n$ follows a shifted and scaled gamma distribution, we obtain, after some manipulation,
\begin{equation*}
    p^*_n(z; \theta) = \frac{\Gamma(n)}{\Gamma\left(\frac{n}{2}\right)^2} \left(\frac{\sqrt{n}\theta^2}{2\sqrt{n}\theta-2z}\right)^{n-1} \int_0^{t_n} \left[x\left(\frac{2}{\theta} - \frac{2z}{\sqrt{n}\theta^2} -x
 \right)\right]^{n/2-1} \mathrm{d}x.
\end{equation*}

Finally, suppose selection acts on $U_n = \sqrt{n}\bar Y_n + W$, where $W\sim N(0, \sigma^2_W)$ for some pre-determined $\sigma_W > 0$, so that the sample gets selected if $U_n< t_n$ (randomization). Then, 
\begin{equation*}
 p^*_n(z; \theta) = \Phi\left\{ \frac{1}{\sigma_W} \left( t_n - \frac{\sqrt{n}}{\theta} + \frac{z}{\theta^2} \right)\right\}.
\end{equation*}
If, instead, we condition on the observed value of the randomized statistic $U_n$, $u_n$ (perhaps due to unavailability of knowledge of the precise threshold $t_n$ applied), the selection function becomes
\begin{equation*}
 p^*_n(z; \theta) = \frac{1}{\sigma_W} \phi\left\{ \frac{1}{\sigma_W} \left( u_n - \frac{\sqrt{n}}{\theta} + \frac{z}{\theta^2} \right)\right\}.
\end{equation*}

\begin{figure}[h]
    \centering
    \includegraphics[width=0.9\linewidth]{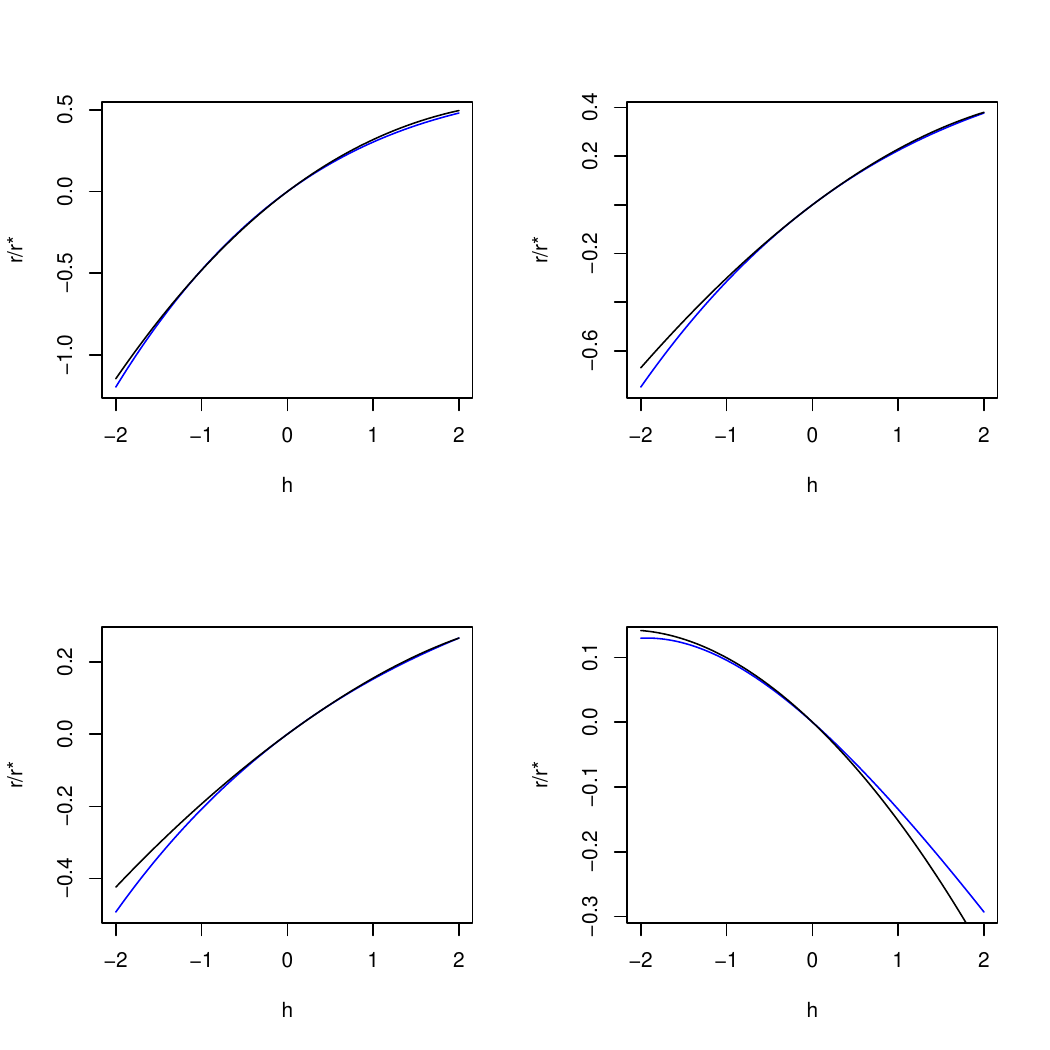}
    \caption{Blue: $r_n(h; \theta)$; black: $r_n^*(h; \theta)$. Top left: deterministic selection with $t_n = 0.5$; top right: data carving with $t_n = 0.5$; bottom left: randomization with $W\sim N(0, 1)$ and $t_n = \sqrt{n}0.5$; bottom right: conditioning on $\sqrt{n}\bar Y_n + W = u_n$, with $W\sim N(0, 1)$ and $u_n = 3.8$.}
    \label{FIG: exponential}
\end{figure}

Figure \ref{FIG: exponential} shows the functions $r_n(h; \theta) = \log\{\varphi_n(\theta + h/\sqrt{n})/\varphi_n(\theta)\}$ and $r^*_n(h; \theta) = \log\{\varphi^*_n(h)/\varphi^*_n(0)\}$ in the four situations considered above, for $\theta = 2$, $n = 40$, and $h\in [-2, 2]$. The main step in the proof of Theorem \ref{SLAN} is to establish that the two functions approximate each other uniformly in compact sets of $h$ as $n\to\infty$. This figure provides visual evidence of the approximation across the range of values of $h$ considered. In all cases, the functions $\varphi^*_n(h)$ where computed by numerical integration of $I_\theta^{1/2}\phi(I_\theta^{1/2}z) p^*_n(z; \theta)$ over $z$. We have seen in Figure 1 how accuracy of the approximation translates to closeness of inferences derived in the true exponential and asymptotically equivalent Gaussian selection models.

\subsection{Inference on winners} \label{SEC: gaussian}
In this section we explore a classic example of selection bias, where inference is conducted on the most promising parameters chosen from a large set of candidates. This constitutes a paradigmatic application of selective inference in the IID setting. In such problems, $N$ independent populations under investigation, with respective parameters $\theta_1, \ldots, \theta_N$, are ranked based on a data-dependent criterion (typically some measure of significance of the parameter), and only the top-ranked populations are selected for further analysis. For example, we may compute $N$ $p$-values of significance and select for formal inference only those parameters whose respective $p$-values are among the smallest $K<N$ or fall below a fixed significance threshold $t$. Once a subset of parameters has been selected, further data may be collected on the selected populations, leading to a data-carving scenario, as described in the previous section. This type of situation is common in clinical studies, where multiple treatments are tested and only the most promising ones are advanced for further investigation. Another common application occurs in genomics, where a vast number of gene expressions are measured, and only those showing the strongest signals, such as significant associations with a disease, are formally analyzed.

We will consider the following specification of the problem. For $j \in \{1, \ldots, N\}$ and $i\in \{1, \ldots, n_1\}$, let $Y_i^j\sim N(\mu_j, \sigma^2_j)$ independently, where all $\mu_j$'s and $\sigma^2_j$'s are unknown. That is, we have random samples of size $n_1$ from $N$ different Gaussian populations. Suppose that we are only interested in means with large values. Accordingly, we compute $N$ $t$-statistics, $\{S_j = \bar Y^j/\sqrt{V^j/n_1} \colon j = 1, \ldots, N\}$, where $\bar Y^j$ and $V^j$ are the maximum likelihood estimators of $\mu_j$ and $\sigma_j^2$, and keep the populations with the largest statistics for inference. We will not analyze any specific choice of selection criterion because, as we shall see, all of them are operationally identical under the considered form of inference. Thus, without loss of generality, suppose that $\mu_j$ is selected for inference if $S_j > T_{n_1, N}(S_{-j}) \equiv T_j$, where $S_{-j}$ contains all $t$-statistics except $S_j$. This threshold could be $T_j = \max(S_{-j})$, or it could be a constant significance threshold obtained from a formal hypothesis test, possibly with a multiplicity correction if $N$ is large. Furthermore, assume that, from each of the selected populations, we collect $n_2$ extra samples to increase inferential power: $Y_{n_1 + 1}^j, \ldots, Y_{n}^j$, where $n = n_1 + n_2$ is the total sample size. This is not necessary for the discussion, but mimics standard procedure in many applications.

In this type of problems, it is common to conduct inference on each of the selected means conditionally on the data from the other $N-1$ populations, in addition to conditioning on the selection event, so as to eliminate the $2(N-1)$ nuisance parameters. For example, let us assume that $\mu_1$ is one of the selected means. Selective inference on $\mu_1$ is then to be conducted from the model $Y_1^1, \ldots, Y_n^1 \sim N(\mu_1, \sigma_1^2)$ with selection event $S_1 > t_1 = T_1(s_2, \ldots, s_N)$, which only depends on the remaining populations through the (now constant) threshold $t_1$. In a selective many-parameter problem, such reduction is often possible and enables analysis of the high-dimensional model in terms of a few independent problems involving low-dimensional selected parameters. Therefore, in our empirical investigations, we shall simply consider selection models of the form $Y_1, \ldots, Y_n \sim N(\mu, \sigma^2)$ with selection event $S > t$, for some fixed $t \in \mathbb{R}$, where the selected population index has been dropped. 

Let $(\bar Y_n, V_n)$ be the maximum likelihood estimator of $\theta = (\mu, \sigma^2)$. We have 
\begin{equation*}
Z_n = (Z_n^1, Z_n^2) \rightsquigarrow N\left(\begin{pmatrix}
0 \\
0
\end{pmatrix}, \begin{pmatrix}
    \sigma^2 & 0 \\
    0 & 2 \sigma^4
\end{pmatrix} \right), 
\end{equation*}
where $Z_n^1 = \sqrt{n}(\bar Y_n - \mu)$ and $Z_n^2 = \sqrt{n}(V_n - \sigma^2 + (\bar Y_n - \mu)^2 )$. The selection function is $p_n(y^n) = \mathbf{1}(s > t)$, where $s$ is the $t$-value as defined above (using only the first $n_1$ observations). 

\begin{figure}[h]
    \centering
    \includegraphics[width=0.9\linewidth]{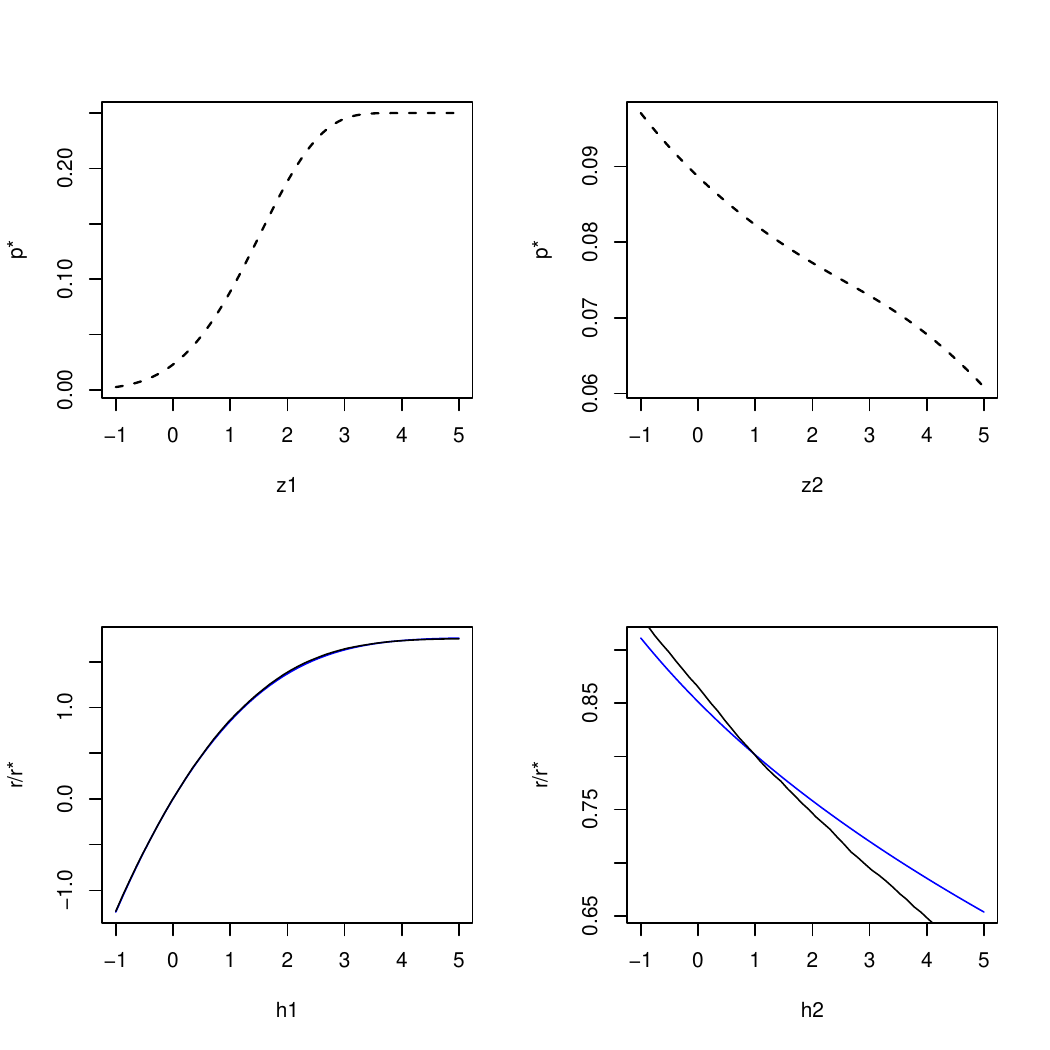}
    \caption{Top left: $p^*_n(z_1, 0; \theta)$; top right: $p^*_n(1, z_2; \theta)$; bottom left: $r_n(h_1, 0; \theta)$ (blue) and $r^*_n(h_1, 0; \theta)$ (black); bottom right: $r_n(1, h_2; \theta)$ (blue) and $r^*_n(1, h_2; \theta)$ (black), with $\theta=(0,1)$.}
    \label{FIG: gaussian}
\end{figure}

If $n_2 = 0$, so that no further data is collected following selection, the asymptotic model for inference admits a tractable analysis. Upon rewriting the MLEs as a function of $(Z_n^1, Z_n^2)$, we obtain the asymptotic selection function
\begin{equation*}
    p_n^*(z; \theta) = \mathbf{1}\left\{ \sqrt{n}\mu + z_1 > t\sqrt{\frac{z_2}{\sqrt{n}} + \sigma^2 - \frac{z_1^2}{n}} \right\},
\end{equation*}
which determines a selection event in the $(z_1, z_2)$-plane. If only $\mu$ is of interest, we can obtain the asymptotic marginal distribution of $Z_n^1$ by integrating out the variance component of the score: 
\begin{equation*}\label{EQ: n2=0}
f(z_1; \theta) \propto \phi\left(\frac{z_1}{\sigma} \right) \left[ \Phi\left\{ U(z_1; \mu, \sigma^2) \right\} - \Phi\left\{ L(z_1; \mu, \sigma^2) \right\} \right],
\end{equation*}
where
\begin{equation*}
L(z_1;\mu, \sigma^2) = \frac{z_1^2}{\sqrt 2\sigma^2 \sqrt{n}} - \frac{\sqrt{n}}{2},
\quad
U(z_1; \mu, \sigma^2) = \frac{\sqrt{n}}{\sqrt 2\sigma^2}\left[\left(\frac{\sqrt{n}\mu + z_1}{t}\right)^2 - \sigma^2 + \frac{z_1^2}{n}\right].
\end{equation*}
If $n_2 > 0$, the asymptotic selection function, $p_n^*(z; \theta) = E_\theta[p_n(Y^n)\mid Z_n = z]$, does not admit, to the best of our knowledge, a simple analytical expression. It is the conditional acceptance probability of a subsample $t$-test given the full-sample sufficient statistics $(\bar Y_n, V_n)$. However, one may compute the asymptotic marginal density numerically:
\begin{equation*}
f(z_1; \theta) \propto \phi\left(\frac{z_1}{\sigma} \right) p_n^*(z_1; \theta), \quad p_n^*(z_1; \theta) = \int_{-\infty}^\infty \phi\left( \frac{z_2}{2 \sigma^2} \right) p^*_n(z_1, z_2; \theta) \mathrm{d}z_2.
\end{equation*}
Figure \ref{FIG: gaussian} shows one-dimensional cuts of $p^*_n(z; \theta)$, $r_n(h; \theta) = \log\{\varphi_n(\theta + h/\sqrt{n})/\varphi_n(\theta)\}$ and  $r_n^*(h; \theta) = \log\{\varphi^*_n(h)/\varphi^*_n(0)\}$ for $n_1 = 20$, $n_2 = 20$, $\theta = (0, 1)$ and $t = 1$. Under the true parameter, the selection probability is $\varphi_n(0, 1) \approx 0.16$. The first two functions were computed via numerical integration, and $r_n^*(h; \theta)$ via Monte Carlo integration of $p^*_n(z; \theta)$. As expected, the approximation of the true selection probability by a Gaussian probability is much more accurate in the direction of the mean ($h_1$) than in the direction of the variance ($h_2$): in this problem $Z_n^1$ is exactly normally distributed. Additionally, since the selection condition is employed to determine large values of $\mu$, the selection function fluctuates much more in the first coordinate, corresponding to $\mu$ (top left plot), than in the second one (top right).

\begin{figure}
    \centering
    \includegraphics[width=\linewidth]{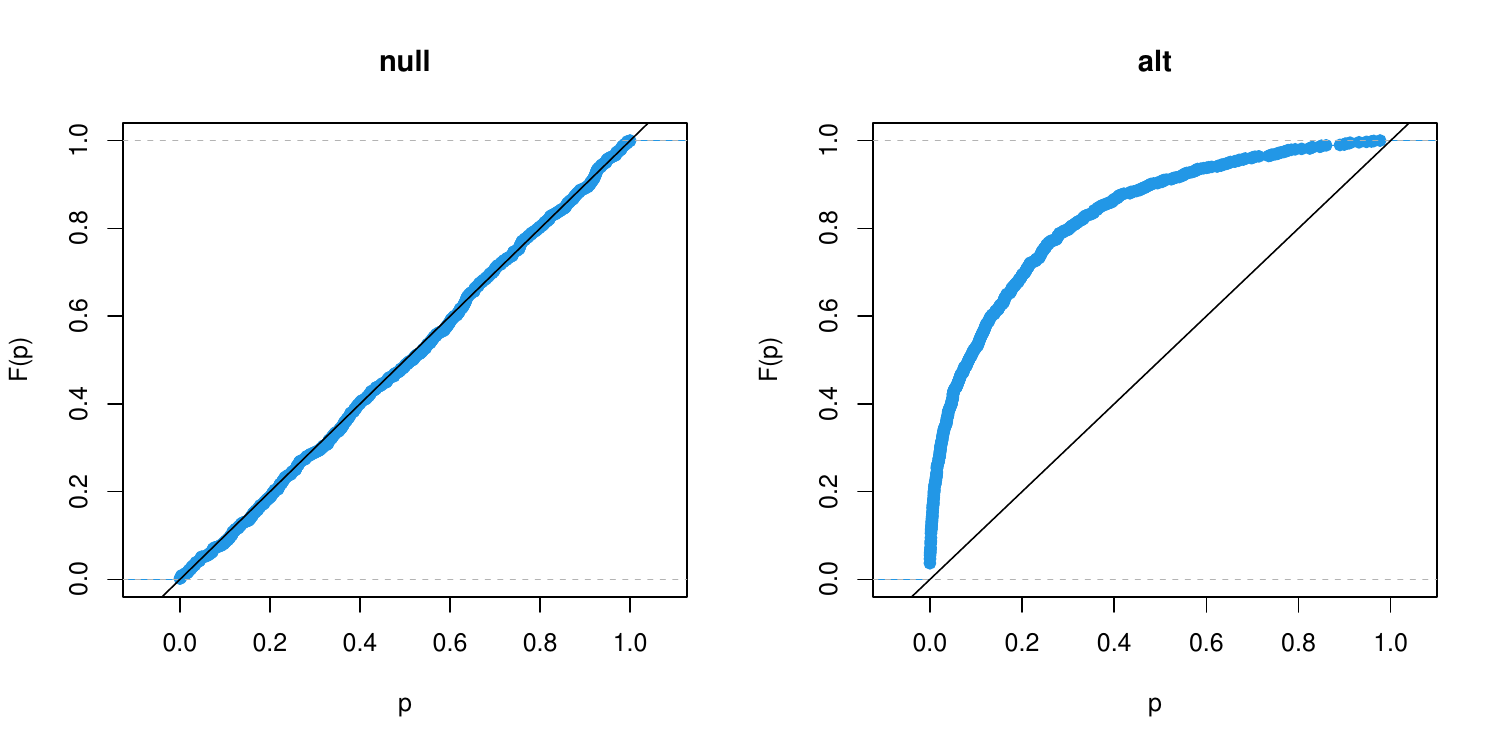}
    \caption{Repeated sampling distribution of $p$-values for the \textit{inference on winners} model, under the null hypothesis (left), and an alternative hypothesis (right); $n_1 = 100$, $n_2 = 50$.}
    \label{Gaussian pvals}
\end{figure}

In both cases, approximate selective $p$-values for $H_0\colon \mu \leq \mu_0$ vs. $H_1\colon \mu > \mu_0$ can be constructed from the asymptotic marginal distribution of $Z_n^1$ at $\mu = \mu_0$, where $\sigma^2$ is substituted by a consistent estimator: 
\begin{equation}\label{EQ: pval}
    p_\text{val}(y^n) = \int_{z_1}^\infty f(z_1; \mu_0, \hat\sigma^2) \mathrm{d}z_1 , \quad z_1 = \sqrt{n}(\bar y_n - \mu_0). 
\end{equation}
 Figure \ref{Gaussian pvals} shows the distribution of \eqref{EQ: pval}, estimated over 1,000 replications, for the case $n _1 = 100$, $n_2 = 50$ and test $H_0\colon \mu \leq 0$ vs. $H_1\colon \mu > 0$. The first plot corresponds to the true generating value $\mu = 0$, and the second one to $\mu = 0.15$. We fixed $\sigma^2 = 1$ and $t = 1$ in both cases. To produce the numerical results we have used the plug-in estimator
\begin{equation*}
    \hat\sigma^2 = \frac{1}{n-1}\sum_{i = 1}^n(Y_i - \bar Y_n)^2,
\end{equation*}
which is consistent as long as the selection probability does not vanish.

This example is considered further in Section 4.2, where it is seen that repeated sampling properties of Bayesian inference under the specified selection model can be well approximated by those in the asymptotically equivalent selective Gaussian location model.

\section{Asymptotic posterior distributions} \label{SEC: Bayes}

Equipped with the formal expansion provided by Theorem \ref{SLAN}, we are now in a position to study the asymptotic shape of posterior distributions constructed from selection models. To this end, suppose we have a prior distribution $\pi(\theta)$ over $\Theta$. Upon observing a sample $Y^n = y^n$ from a selection model $\mathcal{M}_n(u_n)$, the posterior density is given by
\begin{equation} \label{EQ: sel_post}
    \pi_n(\theta\mid y^n) = \frac{\pi(\theta)}{c_n(y^n)} f_\theta(y^n\mid u_n) \propto \frac{\pi(\theta)}{\varphi_n(\theta)} \prod_{i = 1}^n f_\theta(y_i), 
\end{equation}
where $c_n(y^n)$ is the marginal density of $Y^n$ in the selection model. We deliberately avoid the notation $\pi_n(\theta\mid y^n, u_n)$, as it leads to confusion with the standard, non-selective posterior
\begin{equation} 
    \pi_n(\theta\mid y^n, u_n) \propto \pi(\theta \mid u_n) f_\theta(y^n\mid u_n) \propto \pi(\theta ) f_\theta(y^n). 
\end{equation}
We elaborate further on the distinction between both constructions below.

Throughout this section, we adopt the frequentist stance that there is a true underlying data-generating parameter, which is denoted $\theta_0$ so as to avoid confusion with other values of $\theta$ at which the posterior density is evaluated.
Additionally, for a true parameter value $\theta_0 \in \Theta$, we denote by $\pi_{h, n}(h\mid y^n)$ the posterior density for the local parameter $h = \sqrt{n}(\theta - \theta_0)$. The corresponding posterior distribution is
\begin{equation} 
    \Pi_{h, n}(B\mid y^n) = \int_B \pi_{h, n}(h\mid y^n) \mathrm{d}h, \quad B\in \mathcal{B}(\mathbb{R}^p) .
\end{equation}

Posterior distributions of this form are sometimes called \textit{selective posteriors} or \textit{selection-adjusted posteriors}. By contrast to the regular, non-selective IID case, these distributions have a factor of $1/\varphi_n(\theta)$ in their density which, as we shall see, can significantly modify their behavior. Selection-adjusted posterior distributions have appeared in multiple contexts; see \cite{handbook} for a detailed discussion. In selective inference, they arise under the so-called `fixed parameter' regimes, as introduced by \cite{MR1, MR2} and further developed by \cite{yekutieli}. These are notional joint sampling regimes of $(\theta, Y)$ where $\theta$ is sampled from the prior, fixed, and data is generated conditionally on $\theta$ until the selection condition ($U_n = u_n$) is satisfied. Under such sampling regimes it is argued that the posterior density ought to be constructed by combining the prior $\pi(\theta)$ with the likelihood of the conditional model $\mathcal{M}_n(u_n)$, yielding \eqref{EQ: sel_post}. By contrast, `random parameter' regimes arise when $(\theta, Y)$ are sampled jointly until selection occurs. In the latter case, the usual posterior density $\pi_n(\theta\mid y^n)\propto \pi(\theta) \prod_{i = 1}^n f_\theta(y_i)$, blind to selection, is recommended instead for inference \citep{yekutieli}. An interesting conceptual discussion is also provided by \cite{harville}, who notes that \eqref{EQ: sel_post} follows from a standard Bayesian updating of the modified prior $\pi^*(\theta) = \pi(\theta)/\varphi(\theta)$. This has a natural interpretation in the selective inference context: Bayesian inference based on a selection-unadjusted prior has very poor frequentist behavior for $\theta$'s where $\varphi(\theta)$ is small. Thus, by increasing the prior density at those values, this selection-corrected prior provides some level of protection against selection effects. 

Selection-adjusted posteriors are rapidly gaining popularity in theory and practice. Within the regression framework, their utility is highlighted in \cite{panigrahi-scalable}, \cite{panigrahi-grouplasso} and \cite{panigrahi-integrative}. Objective Bayes procedures based on this type of posterior construction are discussed in \cite{woody}. Recent real-data applications include \cite{mackinnon}, who employed them to analyze hop market prices, and \cite{radiogenomics}, who employed them to uncover important gene pathways in a radiogenomic analysis. 

The growing popularity of selection-adjusted posteriors underscores the importance of understanding their theoretical properties, particularly, as we shall see, because their behaviour deviates from that of standard Bayesian models. This section makes two primary contributions towards such understanding. We prove that the consequences of Theorem \ref{SLAN} extend to the Bayesian framework, leading to a result related to the Bernstein von-Mises Theorem for the context of selection. Specifically, we show that the asymptotic shape of the posterior distribution matches that of the corresponding Gaussian selection model under a uniform prior. These latter models provide a more intuitive framework that can be leveraged to better understand the former. Furthermore, as a consequence of the previous result, we indicate a key practical limitation of Bayesian selective inference: probabilistic claims about the parameter (e.g. that it lies within a certain credible interval with probability 90\%), cannot be given a frequentist interpretation with high accuracy. In order to achieve frequentist-matching properties, more complex prior distributions that depend on the sample size are required.

\subsection{Asymptotic behavior of selective posteriors}

By Theorem \ref{SLAN}, under certain conditions, selection models behave asymptotically as selective Gaussian models with observation $Z \sim N(h, I_{\theta_0}^{-1})$ and selection function $p_n^*(z; \theta_0)$. It is thus natural to expect that under a similar set of assumptions the corresponding posterior distributions match asymptotically. To formalize this equivalence, we define the posterior density under the latter model with the improper uniform prior $\pi(h)\propto 1$ as
\begin{equation} \label{EQ: sel_post_norm}
    \pi^*_n(h\mid z) = \frac{\text{det}(I_{\theta_0} )^{1/2} }{\varphi^*_n(h)c^*_n(z)} \phi\left\{  
\Vert I_{\theta_0}^{1/2} (h - z)\Vert \right\}, 
\end{equation}
where $c^*_n(z)$ is a normalizing constant, and denote by $\Pi^*_n(\cdot \mid z)$ the corresponding probability distribution. The following is a direct generalization of Theorem 10.1 in \cite{vandervaart}, and the proof follows largely the same steps. First, we show that the approximation holds conditionally on a sequence of balls of decreasing radius around $\theta_0$, which follows from Theorem \ref{SLAN}. Then, we show that the posterior probability outside these balls is asymptotically negligible, which is guaranteed by the lower-bound assumption on $\varphi_n(\theta_0)$.

\begin{proposition} \label{SBvM}
Consider a sequence of selection models and a true parameter $\theta_0 \in \Theta$ satisfying the assumptions of Theorem \ref{SLAN}, and assume that the prior distribution is absolutely continuous, with a density $\pi(\theta)$ which is continuous and positive at $\theta_0$. Furthermore, suppose that there exists a sequence of tests $T_n\colon \mathcal{Y}^n \to [0, 1]$ such that, for all $\varepsilon > 0$, 
\begin{equation}
    E_{\theta_0}(T_n) \to 0 \quad \text{and} \quad \sup_{\Vert \theta - \theta_0 \Vert \geq \varepsilon} E_{\theta}(1 -T_n) \to 0.
\end{equation}
Then, the sequence of posterior distributions satisfies
\begin{equation}
    \left\Vert \Pi_{h, n}(\cdot \mid Y^n) - \Pi^*_n(\cdot \mid Z_n) \right\Vert_{\mathrm{TV}} \xrightarrow{F_{\theta_0}^{u_n}} 0,
\end{equation}
where $Z_n = n^{-1/2} I_{\theta_0}^{-1}  \sum_{i = 1}^n \triangledown l_\theta(Y_i)  \rvert_{\theta_0}$. That is, the Total Variation distance between the true and asymptotic posteriors vanishes asymptotically in probability conditionally on selection.   
\end{proposition}

When the data-generating parameter is $\theta_0$, Theorem \ref{SLAN} indicates that $Z_n$ follows asymptotically a selective Gaussian distribution with selection function $p^*_n(z; \theta_0)$. Hence, Proposition \ref{SBvM} provides a framework for understanding the frequentist behavior of selective Bayesian methods by examination of their asymptotic, and simpler, Gaussian equivalents. It is important to note that the limiting posterior distributions $\Pi^*_n(\cdot \mid Z_n)$ are not Gaussian unless $p^*_n(z; \theta)$, and thus $\varphi_n^*(h)$, are constant, as selection modifies the likelihood through the factor $1/\varphi_n^*(h)$. The impact of this factor is often manifested in the form of higher posterior variance and heavier tails than those observed in the absence of selection, and, in some cases, pronounced skewness. The latter effect is particularly notable when $\varphi_n^*(h)$ is strictly monotone, reflecting a selection mechanism that systematically favors either large or small values of the location parameter. Additionally, the impact of selection on the posterior grows with the amount of information used for selection. Consequently, hard-truncation mechanisms, such as the ones illustrated in Example \ref{E1}, which employ all available data for selection, produce stronger effects than carved or randomized mechanisms (Examples \ref{E2} and \ref{E3}), as can be observed in Figure \ref{Posterior_density}.

\subsection{Exponential model (ctd.)}

\begin{figure}[h]
\includegraphics[width=\textwidth]{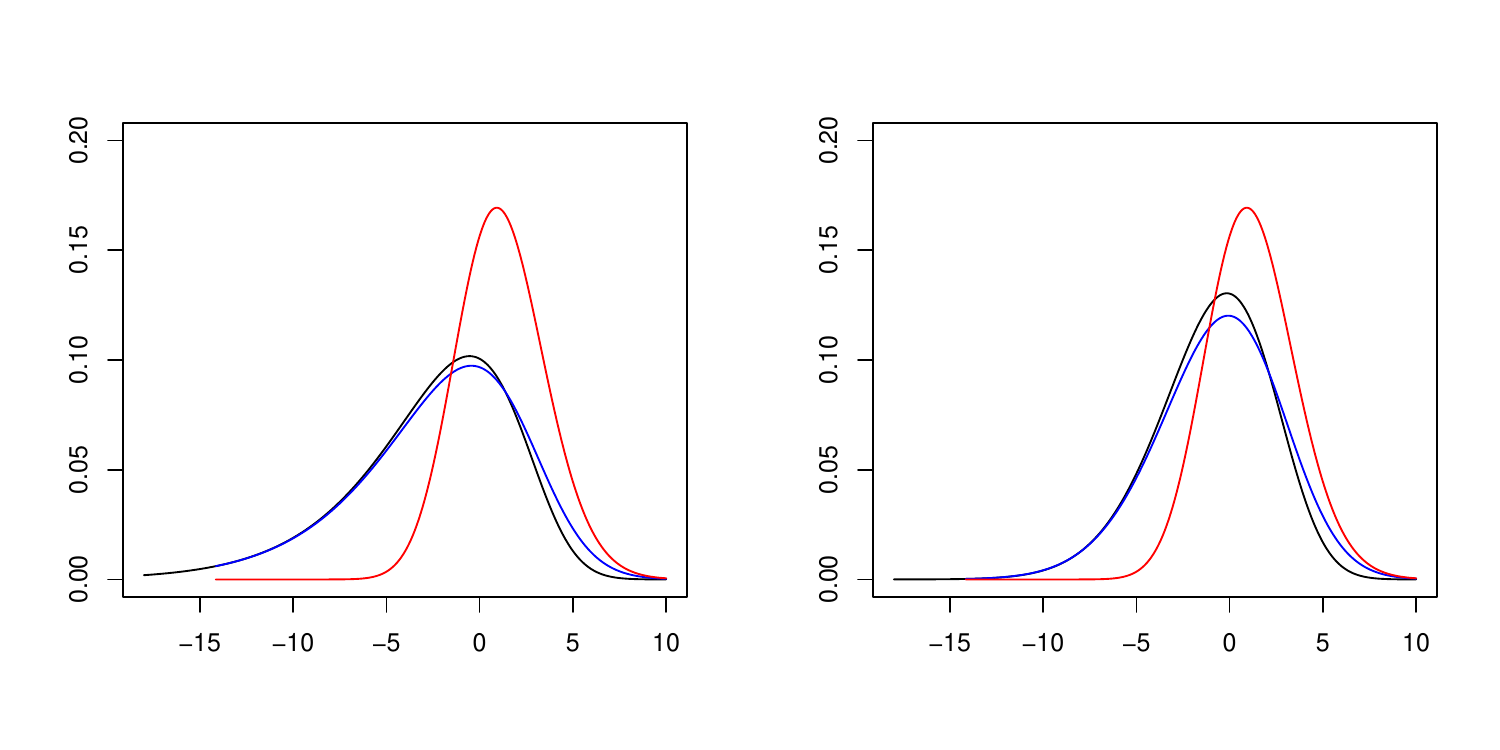}
\caption{Realizations of the posterior densities in the exponential model (blue), their Gaussian approximations (black), and the standard, non-selective posteriors (red). Left: deterministic selection; right, randomized selection.}
\label{Posterior_density}
\end{figure}

Figure \ref{Posterior_density} compares the posterior densities for the local parameter $h = \sqrt{n}(\theta - \theta_0)$, $\pi_{h, n}(h\mid y^n)$, for two exponential models with their theoretical Gaussian limits. We examine two of the scenarios described at the end of Section \ref{SEC: SLAN}: deterministic selection and randomization with Gaussian noise. The model is $Y_1, \ldots, Y_{50}\sim \text{Exp}(\theta)$, $\theta\sim \text{Gamma}(1, 0.1)$, and the true parameter is $\theta_0 = 2$. All posterior densities correspond to an observed sample average of $\bar y_n = 0.45$. The true posterior densities based on the exponential likelihood, $\pi_{h, n}(h\mid y^n)$, are shown in blue, while the Gaussian approximations $\pi^*_{h, n}(h\mid z_n)$, as predicted by the theorem, are shown in black. For comparison, the posterior density in the absence of selection, which is approximately a rescaled $N(1/\bar y_n, \theta_0^2/n)$ under standard theory, is shown in red. Two selection mechanisms are considered: a deterministic one (left), with selection condition $\bar Y_n < 0.5$, and a randomized one (right), with selection condition $\bar Y_n + 0.2 \times N(0, 1 / n) < 0.5$. In both cases, the selection probabilities $\varphi_n(\theta)$ and $\varphi^*_n(h)$ were obtained via numerical integration, as in the numerical investigation presented in Figure \ref{FIG: exponential}. We can readily observe that both posterior densities are non-Gaussian, exhibit higher uncertainty than the non-selective ones, have slowly decaying tails in the direction of the parameter space where $\varphi^*_n(h)$ vanishes, and demonstrate marked skewness, particularly in the case involving deterministic selection.

\subsection{Inference on winners (ctd.)}
\begin{figure}[h]
    \centering
    \includegraphics[width=0.9\linewidth]{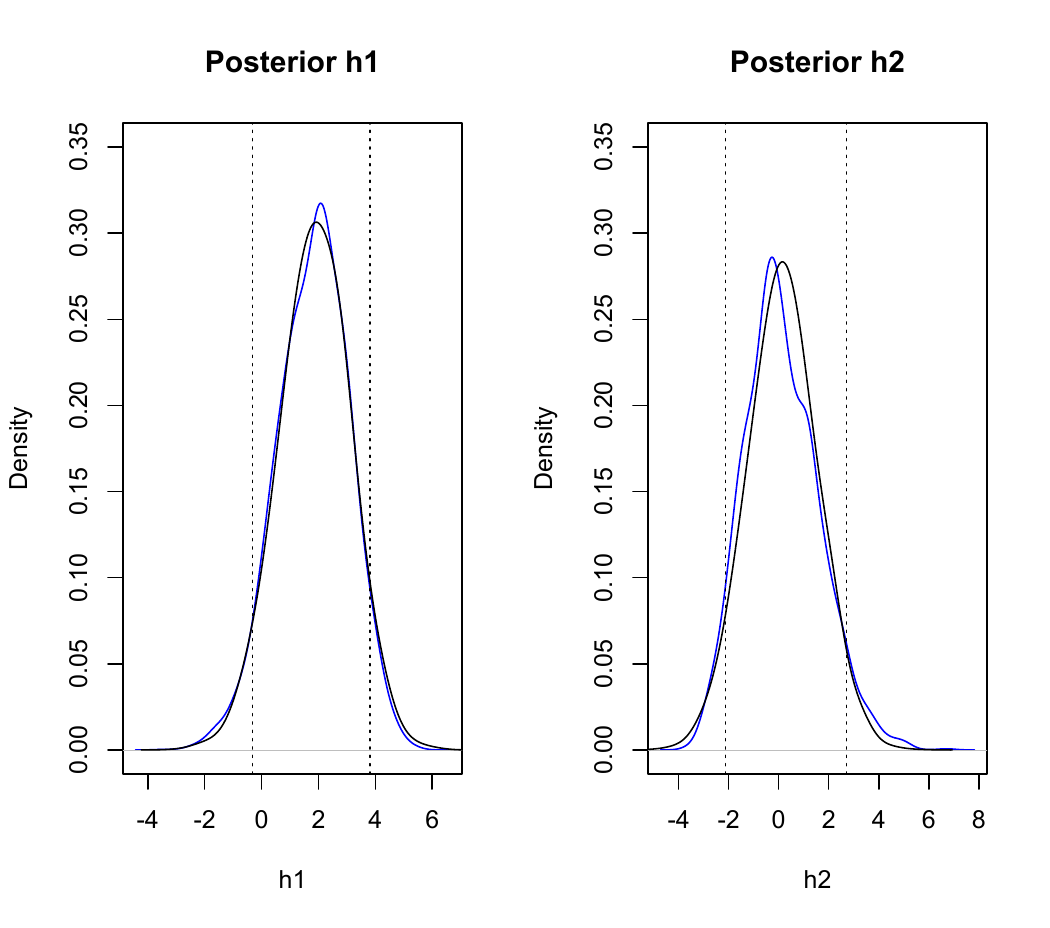}
    \caption{Realizations of the exact marginal posterior densities of $h_1$ and $h_2$ in the inference on winners model (blue) and their Gaussian approximations (black). Exact marginal 90\% credible intervals shown by vertical lines.}
    \label{FIG: Gaussian posteriors}
\end{figure}

For further numerical illustration, we revisit the example discussed in Section 3.5, inference on a normal mean, under a selection condition defined by a $t-$ test on a subsample. We assume true parameter value $\theta_0=(0,1)$ and consider sample sizes $n_1=100, n_2 =50$, with threshold $t_n=1$, which corresponds to a selection probability of about 0.15. Figure \ref{FIG: Gaussian posteriors} illustrates, for a particular data sample, the marginal posteriors for $h_1$ and $h_2$, derived from the exact posterior $\pi_{h, n}(h\mid y^n)$, constructed assuming the prior $\pi(\mu, \sigma^2) \propto 1/\sigma$, and the Gaussian approximation $\pi^*_{h, n}(h\mid z_n)$. In each case, the full posterior for the local parameter $h=(h_1, h_2)$ is constructed by MCMC sampling. The vertical lines in the figure indicate the limits of 90\% credible intervals for $h_1, h_2$ constructed from the exact posterior $\pi_{h, n}(h\mid y^n)$. In the case of $h_1$, for instance, this interval is $(-0.31, 3.81)$ and has probability content under the marginal posterior corresponding to the Gaussian approximation of 0.902. The corresponding figure in the case of $h_2$ is 0.916. The marginal posteriors in the true selection model and the Gaussian approximation model match closely, but not exactly.

We repeated the same analysis for a series of 1000 replications, for several combinations of sample sizes $(n_1, n_2)$ and selection threshold $t_n$. Table \ref{Table: prob content} averages the probability content of the exact 90\% credible interval under the approximate Gaussian posterior over the replications, with the standard deviations over the replications of this given in parentheses. Though the simulation is limited, there are some broad conclusions that can be drawn. The probability content averages are close to 90\%, even for small sample sizes, though there is evidence that a larger sample size is required for matching of the posterior distributions for $h_2$. The probability content figures are more stable over replications for the parameter $h_2$ than the parameter $h_1$, where variability is seen to increase as the selection threshold $t_n$ increases, so that the selection probability shrinks. This is to be expected, as selection primarily affects the mean $\mu$, favouring samples which provide evidence of larger $\mu$, and therefore has greater impact on inference for the local mean $h_1$ than the nuisance parameter $h_2$.

\begin{table}
\caption{Average probability content of marginal 90\% credible intervals, over 1000 replications.}
\label{Table: prob content}
\centering
\begin{tabular}{c c c c c c}
$(n_1,n_2)$ & $t_n$ & \multicolumn{2}{c} {$h_1$} & \multicolumn{2}{c}  {$h_2$} \\ \\
$(20,20)$ &  0 & 0.889&(0.058)& 0.843 &(0.076) \\
        &    1 & 0.895& (0.059)& 0.866 &(0.066) \\
        &    2 & 0.907 &(0.055)& 0.900 &(0.051) \\
        &    3 & 0.918 &(0.064)& 0.926 &(0.069) \\
        \\
$(50,25)$  & 0&  0.897&(0.045)& 0.865& (0.055)\\
          &  1 & 0.910&(0.045)& 0.886 &(0.049) \\
          &  2 & 0.918 &(0.055)& 0.909 &(0.044) \\
          &  3 & 0.920 &(0.086)& 0.919 &(0.054) \\
          \\
$(100,50)$ & 0 & 0.904 &(0.028) & 0.885 &(0.040) \\
        &    1 & 0.911  &(0.028) & 0.890 &(0.038) \\
        &    2 & 0.917 & (0.030) & 0.900 &(0.037) \\
        &    3 & 0.910 &(0.076) & 0.909  &(0.039) \\
\end{tabular}

\end{table}

\subsection{Frequentist mis-calibration of selective Bayesian inference}

In the classical IID framework, the Bernstein-von Mises Theorem establishes asymptotic equivalence between the Bayesian and frequentist approaches to inference on $\theta$. In practical terms, this implies that Bayesian posterior probabilities can be interpreted as frequentist probabilities under the assumption that there is a fixed, true data-generating parameter $\theta_0$. For instance, $1-\alpha$ Bayesian credible sets are also asymptotic $1-\alpha$ frequentist confidence sets. This equivalence stems from the fact that Gaussian location models with a uniform prior on the location satisfy this probability-matching condition exactly.

In the presence of selection, this duality breaks: selection-adjusted posterior distributions can remain systematically mis-calibrated as the sample size increases. This follows by combining Proposition \ref{SBvM} with the observation that the uniform prior on the location parameter lacks a Bayesian-frequentist equivalence under selection, except in the trivial scenarios where the selective distribution is itself Gaussian (Proposition \ref{PROP: miscalibration}). Achieving approximate frequentist calibration in this setting typically requires the use of more carefully tuned priors that account for the selection mechanism and vary with the sample size: see \cite{handbook}. This is not surprising, since, typically in selective inference, the dependence of the selection function on $n$ does not vanish asymptotically. This consideration of mis-calibration is relevant both for Bayesians and for frequentist practitioners who use Bayesian methods for their computational or interpretive convenience.

Formally, a prior density for $\theta$ is (exactly) \textit{probability-matching} if $\Pi(\theta_0\mid Y)\sim U(0, 1)$ under repeated sampling of $Y$ when $\theta_0$ is the true parameter, where $\Pi(\cdot \mid Y)$ is the marginal posterior distribution function of $\theta$ induced by the said prior. This ensures that probabilistic claims about the parameter based on the posterior distribution admit a valid frequentist interpretation.

Conveniently, uniform priors on location parameters are probability-matching. Thus, by the Bernstein-von Mises Theorem, any fixed prior density which is strictly positive at $\theta_0$ is first-order probability-matching, meaning that $\Pi(\theta_0\mid Y) \rightsquigarrow U(0, 1)$ under $\theta_0$, and hence credible sets are asymptotically valid confidence sets. 

By contrast, the following proposition shows that, except in the trivial cases, there is no fixed prior which is probability matching. This analysis encompasses valid improper priors, such as $\pi(\theta)\propto 1$, which arises as the prior density of the asymptotic posterior in Proposition \ref{SBvM}. Consequently, it is often the case that the limiting Gaussian posteriors indicated in Proposition \ref{SBvM} are not calibrated in a frequentist sense. This proposition relates to a well-known one-dimensional result due to \cite{lindley}, which states that the only exact probability-matching priors in one-dimensional models are uniform priors on location parameters. Throughout, for a symmetric matrix $B\in\mathbb{R}^{p\times p}$, we write $B\succ 0$ to mean that $B$ is positive definite, and $B\preceq 0$ to mean that $B$ is negative semidefinite.

\begin{proposition} \label{PROP: miscalibration}
  Let $Y$ follow a selective Gaussian distribution with base model $\{N(\theta, \Sigma)\colon \theta\in \mathbb{R}^p\}$ and selection function $p(y)$, where $\Sigma \succ 0$ is known. A probability-matching prior for $\theta$ exists if and only if the selection function can be written, almost everywhere, as
    \begin{equation*}
        \log p(y) = y^\top A y + b^\top y + c \quad \text{for some } A\in \mathbb{R}^{p\times p}, b\in \mathbb{R}^p, c\in \mathbb{R},
    \end{equation*}
    where $A = A^\top$ and $\Sigma^{-1} - 2 A \succ 0$. In this case, the selective distribution is
    \begin{equation*}
        Y \sim N\left((\Sigma^{-1} - 2 A)^{-1}(\Sigma^{-1}\theta + b), (\Sigma^{-1} - 2 A)^{-1}\right)
    \end{equation*}
    and the probability matching prior is the improper uniform prior on $\theta$, $\pi(\theta)\propto 1$. 
\end{proposition}

If $p$ is additionally required to be bounded, as is the case for the $p^*$ function in Theorem \ref{SLAN}, then the exponential-quadratic representation must satisfy $A \preceq 0$ and $b \perp \ker(A)$, where $\ker(A)=\{v\in\mathbb{R}^p:Av=0\}$. In this case, the selection mechanism can be realized via subspace Gaussian randomization: there exist an integer $r\le p$, a full row-rank matrix $M\in\mathbb{R}^{r\times p}$, a positive definite matrix $\Gamma\in\mathbb{R}^{r\times r}$, and a realization $u\in\mathbb{R}^r$ such that, with
\begin{equation*}
    U = MY + W, \qquad W\sim N(0,\Gamma) \text{ independent of } Y,
\end{equation*}
conditioning on $U=u$ induces the same exponential--quadratic selection function. In this case,
\begin{equation*}
    \log p(y) = -\frac12 y^\top M^\top \Gamma^{-1} M y + (M^\top \Gamma^{-1}u)^\top y  + c \equiv y^\top A y + b^\top y + c,
\end{equation*}
where $A=-(1/2) M^\top \Gamma^{-1}M$ and $b=M^\top \Gamma^{-1}u$ (with $c$ constant given $u$). Additive Gaussian randomization is standard in selective inference; see \cite{splitting}. Cases with $r < p$ correspond to randomization only in certain directions. For example, if $M = (1, 0, \ldots, 0)$, only the first coordinate of $Y$ is randomized.

The degree of departure from uniformity of the marginal posterior density depends on the strength of selection bias, that is, on how much the selection mechanism perturbs the distribution of the data under the true parameter. Miscalibration is typically more pronounced under deterministic selection, and less severe for smooth selection functions, such as those arising from carving and randomized selection mechanisms.

This lack of calibration can be made more precise in the one-dimensional setting when the selection function of the Gaussian model is monotonic. This scenario is of major practical interest as it encompasses selection mechanisms that favor large estimated effects. In such cases, the posterior distribution tends to overstate, under repeated sampling, regions of the parameter space with low selection probability, i.e. it ``overcompensates'' selection bias, as it leads to underestimation, as opposed to overestimation, of the true effect. This result extends Proposition 1 in \cite{handbook}.

\begin{proposition} \label{prop}
Let $Y \sim N(\theta, \sigma^2)$, with $\sigma^2>0$ known, and $p(y)$ a non-constant, increasing, and right-continuous selection function with $\Vert p \Vert_\infty <\infty$. Let $\Pi(\theta\mid  y)$ be the selective posterior distribution function based on the uniform prior $\pi(\theta) \propto 1$. Then,
\begin{equation}
P_{\theta_0}\{\Pi(\theta_0\mid Y) \leq \alpha\mid u\} < \alpha \quad \forall (\alpha, \theta_0) \in (0,1)\times \mathbb{R}.
\end{equation}
\end{proposition}

Uniform miscalibration across the parameter space underscores the inadequacy of the uniform prior for this type of selection. The proof of this result hinges on the existence of a function $\pi(\theta; y)$ which is strictly increasing in $\theta$ for every $y$ such that, if employed as a data-dependent prior density, it achieves exact probability matching. That is, $P_{\theta_0}\{\Pi(\theta_0\mid Y) \leq \alpha\mid u\}= \alpha $ for all possible values of $\theta$ and $\alpha$. Then, by comparison with the constant prior, it follows that the latter produces posterior inferences which are biased towards lower values of $\theta$ under repeated sampling. 

If the selection function is increasing, the posterior distribution tends to place too much probability mass on small values of the parameter relative to what a well-calibrated posterior would. That is, while uncorrected selective inference based on the unconditional model $Y \sim N(\theta, \sigma^2)$ would lead to overestimation of $\theta$, conditioning on the selection event has the opposite effect: it leads to systematic underestimation of $\theta$.

Note that, for a realized $Y = y$, $C(y) = (-\infty, \Pi^{-1}(\alpha\mid y)]$ is a $\alpha 100\%$ credible interval for $\theta_0$. The proposition therefore implies that credible intervals of this form undercover in the frequentist sense, i.e. have coverage strictly smaller than $\alpha$. In \cite{handbook} it is shown that certain default priors, which depend on the sample size and assign low prior probability to parameters with low selection probability, can be used to mitigate this effect. Conversely, the opposite effect occurs in the other direction: the credible intervals $C(y) = [\Pi^{-1}(\alpha\mid y), \infty)$ cover the parameter with probability larger than the nominal level of $1 - \alpha$. In particular, it follows that $p_\text{val}(y) = \Pi(\theta_0\mid y)$ constitutes a valid, though conservative, one-sided $p$-value function for the test $H_0\colon \theta\leq \theta_0$, as $P_{\theta_0}\{p_\text{val}(Y) < \alpha \} < \alpha$.

\begin{figure}[h]
\centering
\includegraphics[width=\textwidth]{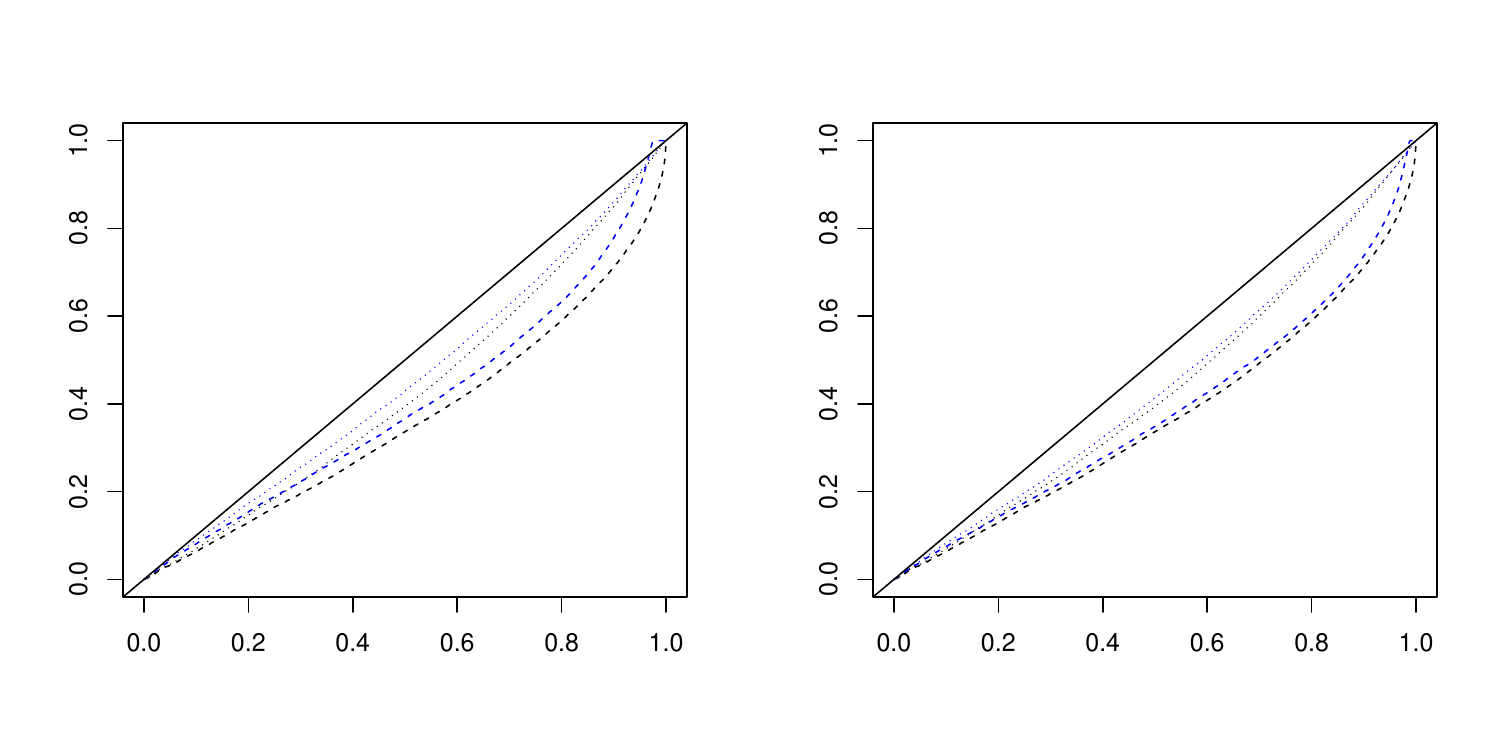}
\caption{Cumulative distribution functions of $\Pi(\theta_0\mid Y^n)$ in the exponential model (blue) with the same settings as before, and of the corresponding Gaussian posterior distributions (black). Left, $n=50$; right, $n=100$. Dashed lines, deterministic selection; dotted lines, randomized selection.}
\label{Posterior_CDF}
\end{figure}

In view of Proposition \ref{SBvM}, a similar behavior ought to arise in more general parametric settings provided the induced selection function on the normalized score, $p^*_n(z; \theta)$, is increasing. Figure \ref{Posterior_CDF} illustrates this phenomenon in the exponential model considered above, where all the induced selection functions are increasing. Qualitatively, we observe that when selection is deterministic there is a large deviation between the distribution of $\Pi(\theta_0\mid Y)$ and the uniform distribution, which would correspond with a perfectly calibrated posterior. This behavior persists from the Gaussian model to the exponential model. Randomizing selection reduces this discrepancy, but the deviation remains significant. 

\section{Discussion}

To date, Gaussianity assumptions have formed the basis of most analysis of inferential properties of selection models and the techniques of selective inference. The asymptotic analysis presented in this paper establishes formal connections between the behavior of selective inference methodologies in more general parametric settings and of their analogues in the corresponding Gaussian limits. An important focus for future work will be to utilize these connections to examine operational and theoretical optimality properties of procedures advocated for inference in the presence of selection. In particular, it will be of methodological importance to explore extensions of the main result to regression settings. 

A key objective of the paper has been to provide insights on the methodological implications of selection for Bayesian inference. Specifically, the convenient asymptotic analogy between frequentist and Bayesian inference in regular random sampling models has been demonstrated to no longer hold under selection. This issue can be especially problematic in settings where the selection mechanism is consistently biased towards one region of the parameter space, as is often the case in practice under common selection conditions. The results presented in this paper point to the need for reexamination of conventional specifications of prior distributions, if appropriate frequentist calibration of the Bayesian inference, as is asymptotically achieved in regular IID problems without selection, is considered important.

Another important avenue for future research concerns the extension of the present results to selection regimes with vanishing selection probabilities. This would provide a unified theoretical basis for the study selection models with IID data, and would substantially broaden the applicability of the framework in high-dimensional contexts.

\section{Acknowledgments}

We would like to thank the Associate Editor and three anonymous referees for their valuable comments and suggestions on a previous version, which have improved the clarity and quality of the work.

\appendix

\newpage

\section{Appendix}

\subsection{Proof of Theorem 1}

Fix a $\theta\in\Theta$ satisfying the required assumptions and let 
\begin{equation*}
 X_n(h) =  \log\{ f_{\theta + h/\sqrt{n}}(Y^n)/f_{\theta}(Y^n)\}.   
\end{equation*}
We have the unconditional expansion
\begin{equation*}
    X_n(h) = h^T I_\theta Z_n - \frac{1}{2} h^T I_{\theta} h + o_{F_\theta}(1),
\end{equation*}
where $Z_n = n^{-1/2} I_\theta^{-1} \sum_{i = 1}^n \triangledown l_\theta(Y_i) \rightsquigarrow N(0, I_\theta^{-1})$. We can write
\begin{equation}
r_n(h) = \frac{\varphi_n(\theta + n^{-1/2} h)}{\varphi_n(\theta)} = E_{\theta}\left[ \exp\left\{  X_n(h) \right\} \mid u_n  \right].
\end{equation}
Consider the sequence of random vectors $[Z, U^*_n]$, where $Z \sim N(h, I_\theta^{-1})$ and $U^*_n\mid Z \sim \text{Bernoulli}\{ p^*_n(Z; \theta)\}$, so that $Z\mid U_n^* = 1$ has selection function $p^*_n(z; \theta)$. The analogous functions for this model are given by
\begin{eqnarray}
     X^*(h) &=& h^T I_\theta  Z - 
\frac{1}{2} h^T I_\theta h; \\
r^*_n(h) &=& \frac{\varphi^*_n(h)}{\varphi^*_n(0)} = E_0\left[ \exp\left\{ X^*(h) \right\} \mid  U^*_n = 1 \right].
\end{eqnarray}
The main part of the proof concerns showing that $r_n(h_n) - r^*_n(h_n) \to 0$ for any convergent $h_n\to h$.
To this end we show the following approximations: 
\begin{eqnarray}
    \text{\textbf{(A)}}& \quad &  E_{\theta}[\exp(h_n^T I_\theta Z_n) p^*_n(Z_n; \theta)] - E_0[\exp(h_n^T I_\theta  Z) p^*_n(Z; \theta)]\to 0; \\
    \text{\textbf{(B)}}& \quad & E_{\theta}[ ( \exp\{  X_n(h_n) \} - \exp\{ h_n^T I_\theta Z_n - (1/2) h_n^T I_\theta h_n \} )p_n(Y^n) ]\to 0. 
\end{eqnarray}

\textbf{(A)} For all $K > 0$, the collection $\{  \exp(h^T I_\theta  Z_n) \colon \Vert h\Vert \leq K, n \in \mathbb{N}\}$ is uniformly integrable under $\theta$. To see this, consider the moment generating function $M(x) = E_{\theta}[ \exp\{ x^T \triangledown l_\theta(Y_1)\}]$. We have $E_{\theta}[ \exp(h^T I_\theta Z_n)] = M(n^{-1/2}h)^n$, so it suffices to prove that 
\begin{equation}
\sup \{  E_{\theta}[ \exp(h^T  Z_n)^2] = M(n^{-1/2} 2 h)^n \colon \Vert h\Vert \leq K, n \in \mathbb{N}\}  < \infty.
\end{equation}
Since $M(0) = 1$ and $\triangledown M(0) = 0$, a Taylor expansion of $\log M(x)$ gives
\begin{equation}
n\log M\left( \frac{2h}{n^{1/2} } \right) = 2 h^T \left[ \triangledown^2 \log M(\tilde x)\right] h, \quad \Vert \tilde{x} \Vert \leq \frac{2 \Vert h \Vert }{n^{1/2}},
\end{equation}
so for some $K' > 0$ and any large enough $n$,
\begin{equation}
\sup_{\Vert h \Vert \leq K, n\in \mathbb{N}} M\left( \frac{2 h}{n^{1/2}} \right)^n \leq \exp \left\{ 2 K^2 \sup_{\Vert x\Vert \leq K'} \Vert \triangledown^2 \log M(x) \Vert_{\text{op}} \right\} < \infty,
\end{equation}
where $\Vert \cdot \Vert_\text{op}$ is the operator norm. Thus, for every $\varepsilon > 0$, there exists a $C_1 > 0$ such that 
\begin{equation}
\sup_{\Vert h\Vert \leq K, n \in \mathbb{N}}E_{\theta}\left[ \exp(h^T  I_\theta Z_n)  \mathbf{1}\{ \exp(h^T I_\theta Z_n) > C_1 \} \right] <  \varepsilon.
\end{equation}
Similarly, there is a $C_2 > 0$ such that
\begin{equation}
\sup_{\Vert h\Vert \leq K} E_0\left[ \exp(h^T I_\theta Z)  \mathbf{1}\{ \exp(h^T I_\theta  Z)  > C_2 \} \right] <  \varepsilon.
\end{equation}
Let $C = \max\{C_1, C_2\}$ and define 
\begin{eqnarray}
g_n(z; h) = \exp(h^TI_\theta z)p^*_n(z; \theta) \mathbf{1}\{ \exp(h^TI_\theta z) \leq C \}, 
\end{eqnarray}
which satisfies $\Vert g_n(\cdot; h)  \Vert_\infty \leq C  \Vert p^*_n(\cdot; h)  \Vert_\infty  $. If the score has an absolutely continuous component, by Theorem 2.5 of \cite{CLT_TV}, $d_{\text{TV}}(Z_n, Z)\to 0$, so $ E_\theta[g_n(Z_n; h)] - E_0[g_n(Z; h)]= o(1)$. Otherwise, define $A(h) =  \{z\colon \exp(h^TI_\theta z) \leq C \}$. We have that
\begin{eqnarray}
\Vert g_n(\cdot; h)_{|A(h) } \Vert_\infty &\leq & C \Vert  p^*_n(\cdot ; \theta)\Vert_\infty ; \\
\Vert g_n(\cdot; h)_{|A(h) } \Vert_\text{L} & \leq & C \left( \Vert  p^*_n(\cdot ; \theta)_{|A(h) }\Vert_\text{L} + \Vert I_\theta h  \Vert \Vert  p^*_n(\cdot ; \theta)_{|A(h) }\Vert_\infty \right) \\
& \leq & C \left( \Vert  p^*_n(\cdot ; \theta)\Vert_\text{L} + \Vert I_\theta h  \Vert \Vert  p^*_n(\cdot ; \theta)\Vert_\infty \right) ;
\end{eqnarray}
where we have used that $\Vert fg\Vert_\text{L} \leq \Vert f\Vert_\infty \Vert g\Vert_\text{L} + \Vert g\Vert_\infty \Vert f\Vert_\text{L}$. Since the bounded-Lipschitz distance metrizes convergence in distribution, we also get $ E_\theta[g_n(Z_n; h)] - E_0[g_n(Z; h)] = o(1)$. Thus, for all $\varepsilon > 0$,
\begin{equation}
\left\vert E_{\theta}[\exp(h_n I_\theta Z_n) p^*_n(Z_n; \theta)] - E_0[\exp(h_n I_\theta  Z) p^*_n(Z; \theta)]\right\vert  \leq 2\varepsilon \Vert  p^*_n(\cdot ; \theta)\Vert_\infty + o(1),
\end{equation}
so the term on the left hand side is $o(1)$, verifying \textbf{(A)}.

In particular, for $h_n = 0$, we get $\varphi_n(\theta) - \varphi_n^*(0) = o(1)$, which implies that $\varphi_n(\theta) / \varphi_n^*(0) = 1 + o(1)$ by the lower bound assumption on $\varphi_n(\theta) $.

\textbf{(B)} First, we have $E_{\theta}[  \exp\{  X_n(h) \}  ] = 1$ for all $h$. Since $Z_n  \rightsquigarrow Z$ and $\{  \exp(h^T I_\theta Z_n) \colon \Vert h\Vert \leq K, n \in \mathbb{N}\}$ is uniformly integrable, it follows that 
\begin{equation*}
    E_{\theta}[  \exp\{ h_n^T I_\theta Z_n - (1/2) h_n^T I_\theta h_n \} ] \to E_0[  \exp\{ h^T I_\theta Z - (1/2) h^T I_\theta h \} ] = 1,
\end{equation*}
so 
\begin{equation}
E_{\theta}[ \exp\{  X_n(h_n) \} - \exp\{ h_n^T I_\theta Z_n - (1/2) h_n^T I_\theta h_n \}  ]\to 0.
\end{equation}
Since $p_n(y^n)$ is bounded, this implies the asserted claim.

Putting all together, using $\dot{=}$ to indicate equality up to an $o(1)$ term, 
\begin{eqnarray}
    &&\varphi_n(\theta)\{ r_n(h_n) - r_n^*(h_n) \} \\
    &=&  E_{\theta}[\exp\{ X_n(h_n)\} p_n(Y^n)] - \frac{\varphi_n(\theta)}{ \varphi^*_n(0)} E_{0}[\exp\{ X^*(h_n) \}  p^*_n(Z; \theta) ]  \\
    &\dot{=}& E_{\theta}[ \exp\{ h_n^T I_\theta Z_n - (1/2) h_n^T I_\theta h_n \} p_n(Y^n) ] - E_{0}[\exp\{ X^*(h_n) \}  p^*_n(Z; \theta) ] \\
    &=& E_{\theta}[ \exp\{ h_n^T I_\theta Z_n - (1/2) h_n^T I_\theta h_n \} p^*_n(Z_n; \theta) ] - E_{0}[\exp\{ X^*(h_n) \}  p^*_n(Z; \theta) ]\\
  &\dot{=}& E_{\theta}[ \exp\{ h_n^T I_\theta Z_n  \} p^*_n(Z_n; \theta) ] - E_{0}[\exp\{ h_n^T I_\theta Z \}  p^*_n(Z; \theta) ]\\
  &=& o(1).
\end{eqnarray}

Since $\varphi_n(\theta)$ is bounded away from zero, we obtain $ r_n(h_n) - r_n^*(h_n) = o(1)$. This establishes the unconditional expansion
\begin{equation}
     \log \frac{f_{\theta + h_n/\sqrt{n}}(Y^n\mid u_n)}{f_{\theta}(Y^n\mid u_n)} = h^T I_\theta Z_n - \frac{1}{2} h^T I_\theta h + \log\frac{\varphi_n^*(0)}{\varphi_n^*(h)} + o_{F_\theta}(1).
\end{equation}
However, the lower boundedness assumption on $\varphi_n(\theta_0)$ ensures that the remainder term is also $o_{F^{u_n}_\theta}(1)$.

Finally, 
\begin{eqnarray}
&& \varphi_n(\theta)\{ E_\theta[g(Z_n)\mid u_n] - E_0[g(Z)\mid U_n^* = 1] \} \\
&=& E_\theta[g(Z_n) p^*_n(Z_n; \theta)] - \frac{\varphi_n(\theta)}{ \varphi^*_n(0)}  E_0[g(Z) p^*_n(Z; \theta)]
\end{eqnarray}
vanishes uniformly for all bounded functions $g$ if Condition 1 is satisfied, showing the last assertion. Otherwise, the same statement holds for all bounded-Lipschitz functions $g$.

\subsection{Proof of Lemma 2}

Under DQM, the non-selective model is Locally Asymptotically Normal at $\theta$, so $f_\theta$ and $f_{\theta+h_n/\sqrt{n}}$ are contiguous for any bounded sequence $h_n$. By Le Cam's first lemma, if $\varphi_n(\theta) = E_\theta[p_n(Y^n)] = o(1)$, $\varphi_n(\theta + h_n/\sqrt{n}) = E_{\theta + h_n/\sqrt n}[p_n(Y^n)] = o(1)$ for any bounded sequence $h_n$. Let $H_n = n^{1/2}(\hat \theta_n - \theta)$. Then, $\varphi_n(\hat\theta_n) = \varphi_n(\theta + H_n/\sqrt{n})$. By conditioning on events $\{\Vert H_n\Vert \leq M\}$ for an arbitrarily large $M$, it is easy to see that $\varphi_n(\theta)\to 0 \Rightarrow \varphi_n(\hat\theta_n) = o_{F^{u_n}_\theta}(1)$. Similarly, if $\varphi_n(\hat\theta_n) = o_{F^{u_n}_\theta}(1)$, then there exists a bounded sequence $h_n$ for which $\varphi_n(\theta + h_n/\sqrt{n}) = o(1)$, so $\varphi_n(\theta) = 0$.

\subsection{Proof of Proposition 3}

The proof follows the arguments of \cite{vandervaart} (page 141) almost step by step, to which we refer the reader for the technical details. Here, we indicate the elements of the proof that require further justification. Let $\Pi_n$ denote the prior distribution for the local parameter $h = \sqrt{n}(\theta - \theta_0)$. First, we need to show that asymptotically the posterior distributions of $h$ obtained with $\Pi_n$ are equivalent to those obtained with the restriction of $\Pi_n$ to $C_n$, the ball of radius $M_n$ around 0, for any $M_n\to\infty$. Then, we apply the local expansion given by Theorem 1.

For the first part, let $F^{u_n}_{n, h}$ denote the selective distribution of the data under $h$. Since $\varphi_n(\theta_0)$ is bounded away from zero, it follows that $F^{u_n}_{n, h_n} \triangleleft \triangleright F^{u_n}_{n, 0}$ for every bounded $h_n$. To see this, note that 
\begin{eqnarray}
&& P_{\theta_0}(A_n\mid u_n) \to 0 \\ &\iff& E_{\theta_0}[p_n(Y^n) \mathbf{1}_{A_n}(Y^n)] \to 0 \\ &\iff& E_{\theta_0 + h/\sqrt{n}}[p_n(Y^n) \mathbf{1}_{A_n}(Y^n)]\to 0, 
\end{eqnarray}
as $F_{n, h_n} \triangleleft \triangleright F_{n, 0}$. If $\varphi_n(\theta_0 + h/\sqrt{n})\to 0$ then $\varphi_n(\theta_0 ) \to 0$ by contiguity, so the last statement is equivalent to $P_{\theta_0 +h/\sqrt{n}}(A_n\mid u_n) \to 0$. Furthermore, we can extend Lemma 10.3 of \cite{vandervaart} to the sequence of selection models. Under the conditions of Corollary 1, for every $M_n\to\infty$, there exists a sequence of tests $T_n$ and a constant $c>0$ such that, for every sufficiently large $n$ and every $\sqrt{n} \Vert \theta - \theta_0 \Vert \geq M_n$,  
\begin{equation}
E_{\theta_0}(T_n ) \to 0 \quad \text{and} \quad E_{\theta}(1 - T_n) \leq e^{-cn(\Vert \theta - \theta_0 \Vert^2 \wedge 1 )}.
\end{equation}
However, under the assumptions on the sequence of selection functions, this also holds conditionally on $u_n$, i.e.
\begin{equation}
E_{\theta_0}(T_n \mid u_n) \to 0 \quad \text{and} \quad E_{\theta}(1 - T_n \mid u_n) \leq e^{-cn(\Vert \theta - \theta_0 \Vert^2 \wedge 1 )}.
\end{equation}
Indeed, we have
\begin{equation}
E_{\theta_0}(T_n \mid u_n) = \frac{1}{\varphi(\theta_0)}E_{\theta_0}[p_n(Y^n) T_n] \to 0,
\end{equation}
as $p_n$ is uniformly bounded and $\varphi_n(\theta_0)$ is bounded away from zero. Moreover, for some constant $K$,
\begin{equation}
E_{\theta}(1 - T_n \mid u_n) = \frac{1}{\varphi(\theta_0)}E_{\theta_0}[p_n(Y^n) (1 - T_n)] \leq K E_{\theta}(1 - T_n ), 
\end{equation}
so the second condition is also satisfied. This establishes the first claim. 

To conclude the proof, let $C$ be a fixed ball around zero of fixed radius. Denote by $\Pi^C_{h, n}(\cdot \mid Y^n)$ and $\Pi^{*C}_n(\cdot \mid Z_n)$ the respective probability measures restricted to $C$, and use analogous notation for the corresponding densities. Their total variation distance can be bounded as
\begin{align*}
& \frac{1}{2}\left\Vert \Pi^C_{h, n}(\cdot \mid Y^n) - \Pi^{*C}_n(\cdot \mid Z_n) \right\Vert_{\mathrm{TV}} \\
& \leq \int \int \left(   1 - \frac{\pi_n(g)  f_{\theta_0+ g/\sqrt{n}}(Y^n\mid u_n) \pi^{*C}_n(h\mid Z_n)  }{\pi_n(h)  f_{\theta_0+ h/\sqrt{n}}(Y^n\mid u_n)\pi^{*C}_n(g\mid Z_n) }  \right)^+ \pi^{*C}_n(g\mid Z_n)   \pi^C_{h, n}(h \mid Y^n)   \mathrm{d}g \mathrm{d}h.
\end{align*}
By Theorem 1 and the regularity assumption on the prior, the integrand converges to zero in probability under the measure $\lambda_C(\mathrm{d}h) F_{\theta_0}^{u_n}(\mathrm{d}y^n) \lambda_C(\mathrm{d}g)$, where $\lambda_C$ denotes the uniform measure on $C$. By the dominated convergence theorem, the right term of the inequality converges to zero in probability under the selective distribution, which concludes the proof, as this holds for an arbitrary $C > 0$.

\subsection{Proof of Proposition 4}

 If a prior is probability-matching for $\theta$, it is also probability-matching for every subparameter of interest $\theta_i$. Without loss of generality, let the parameter of interest be $\theta_1$, $\pi(\theta)$ be a valid arbitrary prior density, not necessarily proper, and denote by $\Pi(\theta_1 \mid y)$ the corresponding marginal posterior distribution function. The prior is probability-matching for $\theta_1$ if $\Pi(\theta_1 \mid Y) \sim U(0, 1)$ under all true parameters $\theta \in \mathbb{R}^p$.
We are going to assume $\Sigma = I_p$, which is analytically easier, and we will consider the general case at the end. 

First, we show that $\Pi(\theta_1 \mid Y) \sim U(0, 1)$ for all $\theta\in \mathbb{R}^p \iff \Pi(\theta_1\mid Y)\mid Y_{-1} \sim U(0, 1)$ almost surely for all $\theta \in \mathbb{R}^p$. 
To see this, write 
\begin{equation*}
P_\theta\{\Pi(\theta_1 \mid Y) \leq \alpha\}  - \alpha= E_\theta[P_{\theta_1}\{\Pi(\theta_1 \mid Y) \leq \alpha \mid Y_{-1} \} - \alpha] = 0.
\end{equation*}
Since $Y_{-1}$ is complete for $\theta_{-1}$ for every fixed value of $\theta_1$ (it is the sufficient statistic of an exponential family with natural parameter $\theta_{-1}$), it follows that 
\begin{equation*}
P_{\theta_1}\{\Pi(\theta_1 \mid Y) \leq \alpha \mid Y_{-1} \}  = \alpha
\end{equation*}
with probability 1 for all $\theta$.

Now, define $H(\theta_1; y) = P_{\theta_1}(Y_1 \geq y_1\mid y_{-1})$. We are going to show that the previous condition can only hold if $\Pi(\theta_1 \mid Y) = H(\theta_1; Y)$ almost surely for all $\theta_1\in \mathbb{R}$. Clearly, $H(\theta_1; Y) \mid y_{-1} \sim U(0, 1)$ for all values of $\theta$ and for all values of $y_{-1}$ for which the conditional distribution is well-defined, and $H(\theta_1; y)$ is a valid distribution function for $\theta_1$ for every fixed $y$, because the conditional model on $Y_1\mid y_{-1}$ is stochastically increasing in $\theta_1$. Additionally, we know that $\Pi(\theta_1\mid y)$ is strictly decreasing in $y_1$ regardless of the prior, as the posterior density can be expressed in the form
\begin{equation*}
\pi(\theta_1\mid y) = \frac{\phi(y_1 - \theta_1) g(\theta_1; y_{-1})}{h(y)}, \quad g(\theta_1; y_{-1}) = \int_{-\infty}^\infty \frac{\pi(\theta)}{\varphi(\theta)}f_{\theta_{-1}}(y_{-1}) \mathrm{d}\theta_{-1}.
\end{equation*}
Let $l_\alpha(\theta, y_{-1})$ be defined by
\begin{equation*}
P_{\theta_1}\{\Pi(\theta_1\mid Y) \leq \alpha\mid y_{-1}\} = P_{\theta_1}\{ Y_1  \geq l_\alpha(\theta, y_{-1})\mid y_{-1}\} =  H(\theta_1;  l_\alpha(\theta, y_{-1}), y_{-1}) = \alpha.
\end{equation*}
For $\alpha = \Pi(\theta_1 \mid y)$ we obtain, from the last equality, $H(\theta_1; y) = \Pi(\theta_1 \mid y)$.

Differentiating this equality with respect to $\theta_1$ and $y_1$, we get
\begin{eqnarray*}
\frac{\partial^2}{\partial \theta_1 \partial y_1} H(\theta_1; y) &=& \frac{\partial}{\partial y_1} \frac{\phi(y_1 - \theta_1) g(\theta_1; y_{-1})}{h(y)} \\
- \frac{\partial}{\partial \theta_1} \frac{\phi(y_1 - \theta_1) p(y)}{\varphi(\theta_1; y_{-1})} &=&  \frac{\partial}{\partial y_1}  \frac{\phi(y_1 - \theta_1) g(\theta_1; y_{-1})}{h(y)},
\end{eqnarray*}
where $\varphi(\theta_1; y_{-1}) = \int_{\mathbb{R}} \phi(y_1 - \theta_1) p(y) \mathrm{d}y_1$, so
\begin{equation*}
 \frac{ p(y)}{\varphi(\theta_1; y_{-1})} \left[ \theta_1 - y_1 + \frac{\partial}{\partial \theta_1} \log \varphi(\theta_1; y_{-1}) \right] = \  \frac{ g(\theta_1; y_{-1})}{h(y)} \left[ \theta_1 - y_1 -  \frac{\partial}{\partial y_1} \log h(y) \right]. 
\end{equation*}



For a fixed $y_{-1}$, rewrite this equation as
\begin{equation*}
    a(\theta_1)\left[ \theta_1 - y_1 + \frac{\partial}{\partial \theta_1} \log \varphi(\theta_1; y_{-1}) \right] = b(y_1) \left[ \theta_1 - y_1 -  \frac{\partial}{\partial y_1} \log h(y) \right],
\end{equation*}
where $a(\theta_1)^{-1} = \varphi(\theta_1; y_{-1})g(\theta_1; y_{-1})$ and $b(y_1)^{-1} = h(y)p(y)$. Differentiating w.r.t. $\theta_1$ we obtain
\begin{equation*}
    a'(\theta_1)\left[ \theta_1 - y_1 + \frac{\partial}{\partial \theta_1} \log \varphi(\theta_1; y_{-1}) \right] + a(\theta_1)\left[ 1 +  \frac{\partial^2}{\partial \theta_1^2} \log \varphi(\theta_1; y_{-1}) \right] = b(y_1) > 0.
\end{equation*}
The left hand side is linear in $y_1$, with slope $-a'(\theta_1) = b'(y_1)$. This implies that $a(\theta_1)$ is linear, but since it is also strictly positive for all $\theta_1$, the only possibility is that $a'(\theta) = b'(y_1)= 0$ for all $\theta_1$. Thus, 
\begin{equation*}
   1 +  \frac{\partial^2}{\partial \theta_1^2} \log \varphi(\theta_1; y_{-1}) \propto 1 \iff
    \frac{\partial^3}{\partial \theta_1^3} \log \varphi(\theta_1; y_{-1}) = 0
\end{equation*}
for all $\theta_1$ and $y_{-1}$ such that $\varphi(\theta_1; y_{-1}) > 0$. By the invertibility of the Weierstrass transform, this condition can be written directly in terms of the selection function, which needs to be almost everywhere quadratic in $y_1$ of the form
\begin{equation*}
    \log p(y) = a(y_{-1})y_1^2 + b(y_{-1})y_1 + c(y_{-1}).
\end{equation*}
Applying the same argument to any orthonormal transformation $\tilde Y = MY\sim N(\tilde \theta, I_p)$, with selection function $p(M^Ty)$, it follows that $\log p(y)$ is quadratic in every possible direction, so 
\begin{equation*}
        \log p(y) = y^T A y + b^T y + c \quad \text{for some } A\in \mathbb{R}^{p\times p}, b\in \mathbb{R}^p, c\in \mathbb{R}.
    \end{equation*}
To conclude, consider a general covariance matrix $\Sigma$ and define $Z = \Sigma^{-1/2} Y\sim N(\tilde \theta, I_p)$, with selection function $\tilde p(z) = p(\Sigma^{1/2}y)$. Probability-matching is invariant to the parametrization, and for the latter model to admit such a prior,
\begin{equation*}
    \log p( y) = \log \tilde p(\Sigma^{-1/2} z) = y^T \Sigma^{-1/2} A \Sigma^{-1/2} y + b^T \Sigma^{-1/2} + c 
 \end{equation*}   
for some $A\in \mathbb{R}^{p\times p}, b\in \mathbb{R}^p, c\in \mathbb{R}$, so the quadratic form holds for $p(y)$ as well.

Since the selection function depends on $A$ via $y^T A y$, only the symmetric part of $A$ is relevant, so we can assume it is symmetric without loss of generality. The further assumption that $\Sigma^{-1} - 2 A \succ 0$ is required for the selective model,
    \begin{equation*}
        Y \sim N\left((\Sigma^{-1} - 2 A)^{-1}(\Sigma^{-1}\theta + b), (\Sigma^{-1} - 2 A)^{-1}\right),
    \end{equation*}
    to be well-defined.

\subsection{Proof of Proposition 5}

Assume for simplicity that $\sigma^2 = 1$, and assume without loss that $\Vert p\Vert_\infty = 1$ (as the model is invariant to rescaling). As in the proof of Theorem 1, consider the model for $(Y, U)$ where $U\mid Y\sim \text{Bernoulli}\{p(Y)\}$, and let $H(\theta; y) = P_\theta( Y\geq y\mid U = 1)$. Note that $H(\theta; Y) \sim U(0, 1)$ under $\theta$. Rewrite this as
\begin{equation}
    H(\theta; y) = \frac{\varphi(\theta; y)}{\varphi(\theta)},
\end{equation}
where $\varphi(\theta; y) = P_\theta(U = 1,  Y\geq y)$. Furthermore, define the function
\begin{eqnarray}
    \pi(\theta; y) = - \frac{H_\theta(\theta; y)}{H_y(\theta; y)},
\end{eqnarray}
where subscripts denote partial differentiation. Formally, $\pi$ can be thought of a data-dependent probability matching prior, in the sense that, when appended to the likelihood of $Y\mid U = 1$, it produces, by construction, a posterior distribution function equal to $H(\theta; y)$, which satisfies the probability-matching condition, as $H(\theta; Y) \sim U(0, 1)$. The proof boils down to showing that this probability-matching prior is increasing in $\theta$ for every fixed $y$.

The cumulant generating function of $Y\mid U = 1$ is $K(t; \theta) = t^2/2 + \theta t + \log\{\varphi(\theta + t)/\varphi(\theta)\}$, so in particular $\varphi_\theta(\theta)/\varphi(\theta) = E_\theta(Y\mid U = 1) - \theta$. Analogously, $\varphi_\theta(\theta; y)/\varphi(\theta) = E_\theta(Y\mid U = 1, Y\geq y) - \theta$, and $\phi(\theta - y)/\Phi(\theta - y) = E_\theta(Y\mid Y\geq y) - \theta$, where $\Phi$ and $\phi$ are the standard Gaussian CDF and PDF, respectively. Using these facts, a direct calculation gives
\begin{equation}
    \pi(\theta; y) \propto \frac{\varphi(\theta;y)}{\Phi(\theta - y)} \times \frac{E_\theta(Y\mid U = 1, Y\geq y) - E_\theta(Y\mid U = 1)}{ E_\theta(Y\mid Y\geq y) - E_\theta(Y)},
\end{equation}
where the proportionality constant depends on $y$ but not on $\theta$. The first factor is increasing in $\theta$ for every $y$, as 
\begin{equation}
    \frac{\partial}{\partial \theta} \log \frac{\varphi(\theta;y)}{\Phi(\theta - y)} = E_\theta(Y\mid U = 1, Y\geq y) - E_\theta(Y\mid Y\geq y) > 0.
\end{equation}
This holds because the selection function is increasing (and non-constant). 

The second factor is also increasing. Since the selection function is right-continuous and $\Vert p\Vert_\infty = 1$, note that the selection condition can be restated as $Y \geq T$, where $T$ is independent of $Y$ and $P(T\leq t) = p(t)$. Write 
\begin{equation}
    E_\theta(Y\mid Y\geq T) = E[E_\theta(Y\mid Y\geq T, T)] = \theta + E\left[ g(\theta - T) \right],
\end{equation}
where $g(x) = \phi(x)/\Phi(x)$, and, similarly,
\begin{equation}
    E_\theta(Y\mid Y\geq T) = \theta + E\left[ g(\theta - \max\{T, y\}) \right].
\end{equation}
We therefore have that
\begin{eqnarray}
        && \frac{E_\theta(Y\mid U = 1, Y\geq y) - E_\theta(Y\mid U = 1)}{ E_\theta(Y\mid Y\geq y) - E_\theta(Y)} \\ &=& E\left[ \frac{g(\theta - \max\{T, y\}) - g(\theta - T)}{g(\theta - y)}\right]\\
        &=& p(y) E\left[1 - \frac{g(\theta - T)}{g(\theta - y)}\mid T < y\right].
\end{eqnarray}
The term inside the expectation is increasing in $\theta$ for every fixed $T < y$ and every $y$ \citep[Proposition 1]{handbook}, which concludes the proof.

\bibliographystyle{abbrvnat}
\bibliography{references}
\end{document}